\documentclass[11pt]{amsart}

\usepackage{amsmath}
\usepackage{amsxtra}
\usepackage{amscd}
\usepackage{amsthm}
\usepackage{amsfonts}
\usepackage{amssymb}
\usepackage{eucal}
\newtheorem{cor}[subsubsection]{Corollary}
\newtheorem{lem}[subsubsection]{Lemma}
\newtheorem{prop}[subsubsection]{Proposition}

\newtheorem{thm}[subsubsection]{Theorem}
\newtheorem{defn}[subsubsection]{Definition}
\theoremstyle{remark}

\newcommand{\nc}{\newcommand}
\nc{\C}{{\mathcal C}}
\nc{\on}{\operatorname}
\nc{\sn}{\mathsf n}
\nc{\pl}{\partial}
\nc{\CB}{{\mathcal B}}
\nc{\CO}{{\mathcal O}}
\nc{\CK}{{\mathcal K}}
\nc{\Spec}{{\operatorname{Spec}}}
\nc{\CG}{{\mathcal G}}
\nc{\CA}{{\mathcal A}}
\nc{\CU}{{\mathcal U}}
\nc{\supp}{\operatorname{supp}}

\nc{\W}{{\mathbf W}}
\renewcommand{\P}{{\mathbf  P}}
\nc{\D}{{\mathbb D}}
\nc{\St}{\operatorname{St}^{\bullet}}
\nc{\CM}{{\mathcal M}}
\nc{\A}{{\mathcal A}}
\nc\GG{\mathbb G}

\nc{\1}{{\mathbf 1}}
\nc{\BC}{{\mathbb C}}
\nc{\CR}{{\mathcal R}}
\renewcommand{\k}{{\mathbb  k}}
\nc{\BQ}{{\mathbb Q}}
\nc{\M}{{\mathcal M}}
\nc{\CL}{{\mathcal L}}
\nc{\CZ}{{\mathcal Z}}
\nc{\U}{{\mathbf U}}
\nc{\tu}{\overset{\bullet}{\mathfrak u}}
\nc{\B}{{\mathcal B}}
\nc{\TM}{\overset{\bullet}{M}}
\nc{\TL}{\overset{\bullet}{L}}
\nc{\BZ}{{\mathbb Z}}

\nc{\CRhom}{\operatorname{RHom}\bul}
\nc{\Fr}{\operatorname{Fr}}
\nc{\Ad}{\operatorname{Ad}}
\nc{\Res}{\operatorname{Res}}
\nc{\gr}{\operatorname{gr}}
\nc{\tr}{\operatorname{tr}}
\nc{\End}{\operatorname{End}}

\nc{\g}{{\mathfrak g}}
\nc{\hatg}{\hat{\mathfrak g}}

\renewcommand{\b}{{\mathfrak b}}
\nc{\sem}{{$\S_{\g}^{\g_{>0}}}}
\nc{\gl}{{\mathfrak g\mathfrak l}}
\nc{\n}{{\mathfrak n}}
\nc{\si}{{\frac\infty 2}}
\nc{\p}{{\mathfrak p}}
\nc{\h}{{\mathfrak h}}

\nc{\Ind}{\operatorname{Ind}}
\nc{\ch}{\operatorname{ch}}
\nc{\Coind}{\operatorname{Coind}}
\nc{\opp}{{\operatorname{opp}}}
\nc{\Ker}{\operatorname{Ker}}
\nc{\im}{\operatorname{Im}}
\nc{\Coker}{\operatorname{Coker}}
\nc{\dirlim}{\underset{\rightarrow}{\operatorname{lim}}}
\nc{\invlim}{\underset{\leftarrow}{\operatorname{lim}}}
\nc{\Sem}{{{\mathbf S}_{\g}^{\g_{>0}}}}
\nc{\CN}{{\mathcal N}}
\nc{\Ext}{\operatorname{Ext}^{\bullet}}
\nc{\ext}{\operatorname{Ext}}
\nc{\tilW}{\til{W}}
\nc{\Mat}{\mathcal{M}at}
\nc{\CV}{\mathcal{V}}
\nc{\lth}{\ell t}
\nc{\BB}{\mathcal{B}}
\nc{\Tor}{\operatorname{Tor}_{\bullet}}
\nc{\tor}{\operatorname{Tor}}
\nc{\Tors}{\operatorname{Tor}_{\frac \infty 2+\bullet}}
\nc{\Exts}{\operatorname{Ext}^{\frac \infty 2+\bullet}}
\nc{\Hom}{\operatorname{Hom}^{\bullet}}
\nc{\ad}{\operatorname{ad}}
\renewcommand{\hom}{\operatorname{Hom}}
\nc{\Tate}{_{\mathsf{Tate}}}
\renewcommand{\mod}{\operatorname{-mod}}
\nc{\Mod}{\operatorname{Mod}}

\nc{\Barb}{\operatorname{Bar}^{\bullet}}
\nc{\F}{{\overline{\mathbb F}_p}}
\nc{\f}{{{\mathbb F}_p}}
\nc{\Q}{{\mathbb Q}}
\renewcommand{\BB}{{\mathbf B}}
\nc{\LL}{{\mathbf L}}
\nc{\upX}{X^{\uparrow}}
\nc{\upcD}{{\mathcal D}^{\uparrow}}
\nc{\upD}{D^{\uparrow}}
\nc{\dX}{X^{\downarrow}}
\nc{\dcD}{{\mathcal D}^{\downarrow}}
\nc{\dD}{D^{\downarrow}}
\nc{\upC}{{\mathcal C}^{\uparrow}}
\nc{\dC}{{\mathcal C}^{\downarrow}}
\nc{\underA}{\underline{A}}
\nc{\underC}{\underline{\CC}}
\nc{\underB}{\underline{B}}
\nc{\underk}{\underline{\k}}
\nc{\Db}{D^{\bullet}}
\nc{\ten}{{\otimes}}
\nc{\tenb}{{\boxtimes}}
\nc{\tenf}{\overset{\operatorname{!}}\ten}
\nc{\tenl}{\overset{\operatorname{L}}\ten}
\nc{\map}{\longrightarrow}
\nc{\eps}{\varepsilon}
\nc{\bs}{\bigskip\\}
\nc{\ms}{\smallskip\\}

\nc{\tilbar}{\widetilde{\operatorname{Bar}}}
\nc{\tilBarb}{\widetilde{\operatorname{Bar}}^{\bullet}}
\nc{\overr}{\overline{R}}
\nc{\overI}{\overline{I}}
\nc{\overX}{\overline{X}}
\nc{\overh}{\overline{h}}
\nc{\overY}{\overline{Y}}
\nc{\overW}{\overline{W}}
\nc{\linbar}{\overline{\operatorname{Bar}}}

\nc{\til}{\widetilde}
\nc{\oppA}{A^{\sharp}}
\nc{\Lemma}{{\bf Lemma:\ }}
\nc{\Theorem}{{\bf Theorem:}\ }
\nc{\Cor}{{\bf Corollary:}\ }
\nc{\Def}{{\bf Definition:}\ }
\nc{\Prop}{{\bf Proposition:}\ }
\nc{\Con}{{\bf Conjecture:}\ }
\nc{\Rem}{{\bf Remark:}\ }
\nc{\dok}{{\em Proof.}\ }

\nc{\SInd}{\operatorname{S-Ind}}
\nc{\SCoind}{\operatorname{S-Coind}}
\nc{\bul}{^{\bullet}}
\nc{\stand}{C^{\frac\infty2+\bullet}}
\nc{\ssn}{\subsection{}}
\nc{\sssn}{\subsubsection{}}
\nc{\codim}{\on{codim}}
\nc{\hgt}{\operatorname{ht}}
\nc{\sqbinom}{\fracwithdelims[][0pt]}

\address{Independent University of Moscow, Pervomajskaya st. 16-18,
Moscow 105037, Russia}

\author{Sergey Arkhipov}

\title{Algebraic construction of contragradient quasi-Verma modules in positive characteristic.}
\date{}

\begin{document} 
\maketitle
\section{Introduction}
The paper grew out of attempts to fill the gap in the proof of Proposition 3.7.1 (ii) in \cite{Ar0}. 
Let us recall the setting there. Denote by $\U_{\CA}$ the Lusztig version of the quantum group
for the root data $(Y,X,\ldots)$ of the finite type $(I,\cdot)$ defined over the algebra $\CA=
\BZ[v,v^{-1}] $ of Laurent polynomials in the variable $v$. Like in~\cite{L1} 
consider the specialization of $\U_\CA$ in characteristic $p$. 
Namely let $\CA_p'$ be the quotient of $\CA$ by the ideal generated by the
 $p$-th cyclotomic polynomial. Then 
$\CA_p'/(v-1)$ is isomorphic to the finite field ${\mathbb F}_p$. 
Thus the algebraic closure $\F$
becomes a $\CA$-algebra. We set $\U_{\F}:=\U_\CA\ten_{\CA}\F$. Let
$\g$ be the semisimple Lie algebra corresponding to the above root data.
 It is known that
the quotient of the algebra $\U_{\F}$ by certain central elements is isomorphic to 
$\U_{\BZ}(\g)\ten \F$, where $\U_\BZ(\g)$ denotes the Kostant integral form 
for the universal enveloping algebra of 
$\g$. 

Recall that an important step in \cite{Ar0} consisted of constructing a certain complex
of $\U_\CA$-modules $\D B\bul_{\CA}(\lambda)$ for a regular dominant integral weight $\lambda$
called the contragradient quasi-BGG complex over $\U_\CA$. The modules in  the
complex are enumerated by the Weyl group in the standard Bruhat order
and are called the contragradient quasi-Verma modules. This complex appears to be a
quantum analogue of a certain classical object. 

Namely consider the Flag variety $\B_{\F}$ for the group $G_{\F}$ and the standard linear bundle
$\CL(\lambda)$ on it. It is known that 
 $\B_{\F}$ is stratified by $B^+_{\F}$-orbits $\{C_{w,\F}|w\in W\}$ called the Schubert cells.
Here and below $W$ denotes the Weyl group corresponding to the root data.
Consider the global Grothendieck-Cousin complex on $\B_{\F}$ with coefficients in $\CL(\lambda)$
with respect to the stratification of the Flag variety by the Schubert cells. We obtain a complex
$$K_{\F}\bul(\lambda)=\underset{m}\bigoplus\underset{w\in
W,\ell(w)=m}{\bigoplus}H^{m}_{C_{ww_0,\F}}(\B_{\F},\CL(\lambda)). $$
Note that by definition the longest element $w_0$ corresponds to the largest Schubert cell and the unit
element $e\in W$ corresponds to the zero dimensional Schubert cell in $\B_{\F}$.
An easy observation shows that the complex carries a natural action of  $\U_{\F}(\g)$. 

So our aim was to compare the complexes $K_{\F}\bul(\lambda)$ and 
$\D B\bul_{\CA}(\lambda)\ten_{\CA}\F=:\D B\bul_{\F}(\lambda)$. The main problem 
here is that the first complex is constructed in purely algebro-geometric terms while the second
one is constructed algebraically ``by hands''.

Thus the present paper is devoted to spelling out a construction of the complex $K_{\F}\bul(\lambda)$
in purely algebraic terms. One of the two main results   of the 
paper (Theorem~\ref{compare}) provides an 
isomorphism between $H^{\ell(w)}_{C_{ww_0,\F}}(\B_{\F},\CL(\lambda))$ 
and a suitably defined
characteristic $p$ version of the contragradient quasi-Verma module 
$\D\til{M}_{\F}(w\cdot\lambda)$.
We also compare the differentials in the complexes $\D B\bul_{\F}(\lambda)$ and 
$K_{\F}\bul(\lambda)$. The second main result of the paper (Theorem~\ref{main}) claims that 
the differential s in the complexes coincide. This means that the complex $\D B\bul_{\F}(\lambda)$
is isomorphic to the global Grothendieck-Cousin complex on the flag variety with coefficients in
$\CL(\lambda)$.

\bigskip

Let us describe the structure of the paper. 

The second section can be thought of as a preparatory one. There we introduce the main objects of our
 study in characteristic $0$. Namely for any element of the Weyl group $w\in W$ we define a
{\em semiregular bimodule} ${\mathcal S}_w$ over the Lie algebra $\g$. Using this bimodule we 
construct a kernel functor $S_w:\ \g\mod\map\g\mod$. Composing this functor with the 
functor $T_w$ of twisting the $\g$-action by adjunction with a certain lift of $w\in W$ to the group $G$
we obtain the functor $\Theta_w:\ \g\mod\map\g\mod$. It turns out that $\Theta_w$ 
possesses some nice properties. In particular it preserves the standard category ${\mathcal O}(\g)
\subset \g\mod$. When restricted to the bounded derived category of the category  ${\mathcal O}(\g)$
the derived functor of $\Theta_w$ becomes an autoequivalence. 

We define the twisted Verma modules over $\g$ as the images of the usual Verma modules under 
the functor $\Theta_w$. Among the twisted Verma modules we are mostly interested in 
the contragradient {\em quasi-Verma modules} defined as images of the antidominant 
Verma modules under $\Theta_w$. However the picture is quite simple in characteristic $0$: we prove
that the contragradient {\em quasi-Verma modules} coincide with the usual contragradient Verma 
modules.

In the third section we show that the picture in characteristic $p$ is much more complicated.
We start the section with defining  the left semiregular module ${\mathcal S}_{w,\F}$ over the
hyperalgebra $\U_{\F}(\g)$. We show that the $\U_{\F}(\g)$-module ${\mathcal S}_{w,\F}$
is isomorphic to the top local cohomology space $H^{\codim N_w}_{N_{w,\F}}
(G_{\F},{\mathcal O}_{G_{\F}})$. This makes ${\mathcal S}_{w,\F}$ a $\U_{\F}(\g)$-bimodule.

Again, like in characteristic $0$, 
using this bimodule we construct a kernel functor $S_{w,\F}:\ \U_{\F}(\g)\mod\map\U_{\F}(\g)\mod$.
After adding some $w$-twist of the $\U_{\F}(\g)$-action we obtain the twisting functor
 $\Theta_{w,\F}:\ \U_{\F}(\g)\mod\map\U_{\F}(\g)\mod$. We define 
the twisted Verma modules over $\U_{\F}(\g)$ as the images of the usual Verma modules  over
 $\U_{\F}(\g)$ under 
the functor $\Theta_{w,\F}$. Again we  are mostly interested in the contragradient 
{\em quasi-Verma modules} being the images of the antidominant Verma modules 
over $\U_{\F}(\g)$ under  the functor $\Theta_{w,\F}$. For a regular dominant integral weight 
$\lambda$ we denote $\Theta_{w,\F}(M_{\F}(w_0\cdot\lambda) )$ by 
$\D\til{M}_{\F}(ww_0\cdot\lambda)$.

 Our main result in the third section claims that
$\D\til{M}_{\F}(ww_0\cdot\lambda)$ is isomorphic to the local cohomology space of the Flag variety
$\B_{\F}$ with support in the Schubert cell corresponding to $w$ and with coefficients in $\CL(\lambda)$.

This means in particular that the underlying $\U_{\F}(\g)$-module of the complex 
$K_{\F}\bul(\lambda)$
can be written as follows:
$$
K_{\F}\bul(\lambda)=\bigoplus K_{\F}^m(\lambda),\  K_{\F}^m(\lambda)=
\underset{\ell(w)=m}\bigoplus\D\til{M}_{\F}(w\cdot\lambda).
$$

The fourth section is devoted to an algebraic construction of the differentials in $K_{\F}\bul(\lambda)$.
We begin with the easiest case of the map 
$$
K_{\F}^{\dim\B-1}(\lambda)\map K_{\F}^{\dim\B}(\lambda):\ 
\underset{i\in I}\bigoplus\ \D\til{M}_{\F}(s_iw_0\cdot\lambda)
\map\D\til{M}(w_0\cdot\lambda).
$$
(see Theorem~\ref{square1}). We show that the above maps can be constructed starting from the 
differential in the Grothendieck-Cousin complex in the $SL(2)$ case using the induction functor
$\on{Ind}_{\U_{\F}(\p_i)}^{\U_{\F}(\g)}$. Here $\p_i$ denotes the $i$-th standard parabolic 
subalgebra in $\g$.

The next case we treat is the one of a simple reflection. We construct ``by hands'' the component of 
the differential $d_{w}:\ \D\til{M}(w\cdot\lambda)\map\D\til{M}(ws_i\cdot\lambda)$, where $i\in I$ and
$\ell(ws_i)=\ell(w)+1$ (see Theorem~\ref{square2}). We show that the map $d_{w}$ can be obtained 
from the appropriate differential component in the previous case by applying the twisting functor
$\Theta_{w,\F}$.

Finally we treat the general case of a component $\D\til{M}_{\F}(w\cdot\lambda)\map
\D\til{M}_{\F}(w'\cdot\lambda)$, where $\ell(w')=\ell(w)+1$ and $w'$ follows $w$ 
in the Bruhat order on the Weyl group. Considering the contragradient dual complex 
$\D K_{\F}\bul(\lambda)$ we prove essentially by a deformation argument that certain chains of the
differential components define canonical embeddings $\til{M}_{\F}(w\cdot\lambda)\hookrightarrow
\til{M}_{\F}(\lambda)$. Moreover for any pair of elements $w,w'$ of the 
Weyl group such that $w'$ follows $w$ in the Bruhat order
the submodule $\on{Im}\left((\til{M}_{\F}(w'\cdot\lambda)\right)\subset
\til{M}_{\F}(\lambda)$ is embedded
into the submodule $\on{Im}\left(\til{M}_{\F}(w\cdot\lambda)\right)\subset\til{M}_{\F}(\lambda)$
(see Theorem~\ref{embeddings}).

In particular the differential component $\D\til{M}_{\F}(w\cdot\lambda)\map
\D\til{M}_{\F}(w'\cdot\lambda)$ can be characterized implicitly this way.

The fifth section is devoted to comparing the complexes $\D B\bul_{\F}(\lambda)$ and 
 $K_{\F}\bul(\lambda)$. We begin the section with recalling the construction of the complex
$ B\bul_{\CA}(\lambda)$. Then we consider the specialization of the complex  
$ B\bul_{\F}(\lambda)$ over $\F$ and prove the main result of the paper.
\vskip 1mm
\noindent
{\bf Main Theorem:} For a regular dominant integral weight $\lambda$ 
the complexes $\D B\bul_{\F}(\lambda)$
and $K_{\F}\bul(\lambda)$ are isomorphic as complexes of $\U_{\F}(\g)$-modules.
\vskip 2mm
\noindent
{\bf Acknowledgments.} The paper was bourn in numerous discussions with Prof. Henning Haahr 
Andersen and Prof. Niels Lauritzen during the author's visit to Aarhus University, Denmark. I am 
happy to express my deepest gratitude to them for their inspiring interest and generous exchange of 
ideas. Also I would like to thank Aarhus University, Denmark, and Max Planck Institute for 
Mathematics, Bonn, Germany, for hospitality and excellent working conditions.

\section{Characteristic zero case}
\subsection{Semiregular bimodules over a semisimple Lie algebra. Algebraic setting.}
Let $G$ be a semi-simple algebraic group over $\BC$. Let $T$ be the Cartan group of $G$ 
and let $(I,X,Y)$ be the corresponding root data, i.e $X$ is the set of characters
$T\map \GG_m$ (i.e. the weight lattice of $G$) and $Y$ is the set
of cocharacters $\GG_m\map T$ (i.e. the coroot lattice of $G$).
Denote the Weyl group 
(resp. the root system, resp. the positive root system) for $(I,X,Y)$
by $W$ (resp. by $R$, resp. by $R^+$). Let $\rho$ denote the half-sum of the positive roots.
 We use the notation $\ell(\cdot)$ for the usual length function on the Weyl group.

 We denote by $\g=\on{Lie}(G)$
the corresponding semisimple Lie algebra. By definition it is graded by the coroot 
lattice. Choose $w\in W$.  We make $\g$ a graded Lie algebra by setting the grading of all
the simple coroots $\check\alpha_i$ such that $w(\check\alpha_i)$ is a  positive coroot as well
equal to $1$ and putting the grading of the rest of the simple coroots equal to $0$.
 Thus $\g=\underset{m\in\BZ}\bigoplus\g_m^w$. Denote
by $\n_w^-$ the  negative triangular subalgebra,
$\n_w^-=\underset{m<0}\bigoplus\g_m^w$ 
The standard Cartan subalgebra $\g_0=\g_0^e$, where $e\in W$ denotes the unity element,
 is denoted by $\h$.

\begin{defn}
The left graded $\g$-module ${\mathcal S}_{\n_w^-}:=U(\g)\ten_{U(\n_w^-)}U(\n_w^-)^*$, where
by definition $U(\n_w^-)^*=\underset{m\in\BZ}\bigoplus\left(U(\n_w^-)_{-m}\right)^*$, is called
{\em the left semiregular $\g$-module} with respect to the nilpotent subalgebra $\n_w^-$.
\end{defn}

Below we show that the left module ${\mathcal S}_{\n_w^-}$ carries a natural structure of a 
$\g\oplus\g$-module. In fact the existence of this structure follows implicitly from 
a deep result in \cite{Ar2}. We prefer not to use the result and to introduce
the structure explicitly.  In fact we work in a more general setup 
including in particular the one stated above.

Let $\g=\underset{m\in\BZ}{\bigoplus}\g_m$ be a graded Lie
algebra with a finite dimensional negatively graded Lie subalgebra
$\n\subset\g$ that acts locally ad-nilpotently both on $\g$ and on $U(\g)$.

\subsubsection{The abelian case.}
Suppose in addition that $\n$ is abelian and one dimensional, $\n=\BC e$. Then $U(\n)=\BC [e]$ and 
it will be convenient for us to identify $\BC [e]^*$ with the $\BC [e]$-module $\BC [e,e^{-1}]/\BC [e]$.
In particular the left semiregular module is isomorphic to 
$$
U(\g)\ten_{\BC [e]}\left(\BC [e,e^{-1}]/
\BC [e]\right)=
U(\g)\ten_{\BC [e]}\BC [e,e^{-1}]/U(\g).
$$
In fact below we show that the $\g$-module $U(\g)_{(e)}=
U(\g)\ten_{\BC [e]}\BC [e,e^{-1}]$ carries a natural structure of an associative algebra. 
Indeed, we have the 
following statement.

\begin{lem}
For any $u\in U(\g)$ the series
$$
[e^{-1},u]=\left[\frac 1{1-(1-e)},u\right]=\sum\ad_{1-e}^m(u)
$$
is finite and provides a well defined element of $U(\g)_{(e)}$.\qed
\end{lem}

\begin{cor}
 $U(\g)_{(e)}$ becomes an associative algebra with the natural subalgebra $U(\g)\subset U(\g)_{(e)}$.
In particular the left semiregular module ${\mathcal S}_{\n}=U(\g)_{(e)}/U(\g)$ carries a natural 
structure of a $\g\oplus\g$-module.
\qed
\end{cor}

\subsubsection{General case.}                \label{filtration}
Let $\g$ and  $\n$ be like in the beginning of this subsection.
Then there exists a filtration $F$ on $\n$:
$$
\n=F^0\n\supset F^1\n\supset\ldots\supset F^{\operatorname{top}}\n=0
$$
such that each $F^i\n$ is an ideal in $F^{i-1}\n$ and there exist abelian one dimensional 
subalgebras $\n^i\subset F^i\n$ such that we have $F^i\n=\n^i\oplus F^{i+1}\n$
as a vector space.

\begin{thm}                                              \label{main2}
There exists  an inclusion of algebras
$U(\g)\hookrightarrow\End_{\g}({\mathcal S}_{\n})$ such that the image of $U(\g)$ is
a dense subalgebra in $\End_{\g}({\mathcal S}_{\n})$.
\end{thm}
\vskip 1mm
\noindent
\Rem
The rest of this subsection is devoted to the proof of the above Theorem.

\begin{proof}
Suppose $\g\supset\n$ and $\n={\n^+}\oplus{\n^-}$ as a graded vector space, where ${\n^-}$ is
an ideal in $\n$ and ${\n^+}$ is a subalgebra in $\n$.  Consider the left
$\n$-modules
$$
{\mathcal S}_{{\n^-}}^{\n}:=U(\n)\ten_{U({\n^-})}U({\n^-})^*\text{ and }
{\mathcal S}_{{\n^+}}^{\n}:=U(\n)\ten_{U({\n^+})}U({\n^+})^*.
$$
Then the $\g$-modules ${\mathcal S}_{{\n^-}}$ and $U(\g)\ten_{U(\n)}{\mathcal S}_{{\n^-}}^{\n}$
(resp.
the $\g$-modules ${\mathcal S}_{{\n^+}}$ and $U(\g)\ten_{U(\n)}{\mathcal S}_{{\n^+}}^{\n}$)
are naturally isomorphic.

The induction functor $\Ind_{{\n^-}}^{\n}:\ {\n^-}\mod\map\n\mod$
takes $U({\n^-})^*$ to ${\mathcal S}_{{\n^-}}^{\n}$ and provides an inclusion of algebras
$U(\n)\subset\End_{\n}({\mathcal S}_{{\n^-}}^{\n}$):
$$
n\cdot a\ten f:=a\ten n\cdot f,\text{ where }n\in{\n^-}, f\in U({\n^-})^*,a\in U(\n),
\ n\cdot f(n'):=f(n'n), \n'\in{\n^-}.
$$
\sssn
Consider the following action of ${\n^+}$ on $U(\n)\ten U({\n^-})^*$:
\begin{gather*}
n^+\cdot a\ten f:=an^+\ten f+a\ten[n^+,f],\text{ where }
n^+\in{\n^+},f\in U({\n^-})^*,a\in U(\n),\\
[n^+,f](u):=f([u,n^+]),\ u\in U({\n^-}).
\end{gather*}
\begin{lem} \label{act}
The  action of ${\n^+}$  on $U(\n)\ten U({\n^-})^*$ introduced above
commutes with the left regular action of $U(\n)$.
Moreover it is well defined on ${\mathcal S}_{{\n^-}}^{\n}$. Along with the action of
${\n^-}$ on ${\mathcal S}_{{\n^-}}^{\n}$ it defines the inclusion of algebras $U(\n)\subset
\End_{\n}({\mathcal S}_{{\n^-}}^{\n})$.
\end{lem}

\begin{proof}
The first statement is obvious. Let $a\in U(\n)$, $f\in U({\n^-})^*$,
$n^+\in{\n^+}$, $n^-\in{\n^-}$. Note that $[n^-,n^+]\in{\n^-}$. Thus we have
\begin{gather*}
n^+\cdot an^-\ten f=
an^-n^+\ten f+an^-\ten[n^+,f]=
an^-n^+\ten f+a\ten n^-[n^+,f]\\=
an^+n^-\ten f+a\ten [n^-,n^+]f-a\ten [n^-,n^+]f+a\ten [n^+,n^-f]=
n^+\cdot a\ten n^-f.
\end{gather*}
To prove the third statement note that
\begin{gather*}
n^+\cdot n^-\cdot a\ten f=
an^+\ten n^-f+a\ten[n^+,n^-f]\\=
an^+\ten n^-f+a\ten n^-[n^+,f]+a\ten[n^+,n^-]f=
n^-\cdot n^+\cdot a\ten f+[n^+,n^-]\cdot f.
\end{gather*}
\end{proof}
The following statement is crucial in the proof of the Theorem.

\begin{lem}
There exists an isomorphism of left $\n$-modules 
$$U(\n)^*\cong {\mathcal S}_{{\n^-}}^{\n}
\ten_{U({\n^+})}U({\n^+})^*.$$
\end{lem}

\begin{proof}
Let $f^+\in U({\n^+})^*$, $f^-\in U({\n^-})^*$. Denote the element
$(1\ten f^-)\ten f^+\in {\mathcal S}_{{\n^-}}^{\n}\ten_{U({\n^+})}U({\n^+})^*$ by $f^-\ten
f^+$.  Then such elements form a base of the vector space
${\mathcal S}_{{\n^-}}^{\n}\ten_{U({\n^+})}U({\n^+})^*\cong U({\n^-})^*\ten U({\n^+})^*$. 
We calculate
the action of $\n$ on ${\mathcal S}_{{\n^-}}^{\n}\ten_{U({\n^+})}U({\n^+})^*$ in this base.  For
$n^-\in{\n^-}$ we have
$$
n^-(f^-\ten f^+)=(n^-f^-)\ten f^+.
$$
For $n^+\in{\n^+}$ we have
\begin{gather*}
n^+(f^-\ten f^+)=(n^+\ten f^-)\ten f^+=(n^+\cdot(1\ten f^-))\ten f^+
-[n^+,f^-]\ten f^+\\ =
f^-\ten n^+f^+-[n^+,f^-]\ten f^+.
\end{gather*}
Now recall that ${\n^-}$ is an ideal in $\n$ and the actions  of the subalgebras
${\n^-}$ and ${\n^+}$ on $U(\n)^*\cong U({\n^-})^*\ten U({\n^+})^*$ are given by these
very formulas.  
\end{proof}

We return to the situation of (\ref{main2}). Note that each module
${\mathcal S}_{\n^k}$ is in fact a $\g$-bimodule by Lemma~\ref{act}.

\begin{lem}  \label{iterate}
The left $\g$-modules ${\mathcal S}_{\n}$ and ${\mathcal S}_{\n^0}\ten_{U(\g)}\ldots\ten_{U(\g)}
{\mathcal S}_{\n^{\operatorname{top}-1}}$ are naturally isomorphic to each other.
\end{lem}

\begin{proof}
We prove by induction by $\operatorname{top}-k$ that
the $\g$-modules ${\mathcal S}_{F^k\n}$ and ${\mathcal S}_{F^{k+1}\n}
\ten_{U(\g)}{\mathcal S}_{\n^k}$
are naturally isomorphic to each other.
Note that each time the induction hypothesis provides the $\g$-bimodule
structure on ${\mathcal S}_{F^{k+1}\n}$:
$$
U(\g)\hookrightarrow\End_{\g}({\mathcal S}_{F^{\operatorname{top}-1}\n})
{\overset{\cdot\ten_{U(\g)}{\mathcal S}_{\n^{\operatorname{top}-2}}}
{\map}}
\ldots
{\overset{\cdot\ten_{U(\g)}{\mathcal S}_{\n^{k+1}}}
{\map}}
\End_{\g}({\mathcal S}_{F^{k+1}\n}).
$$
Thus by the previous Lemma we have
$$
{\mathcal S}_{F^{k+1}\n}\ten_{U(\g)}{\mathcal S}_{\n^k}\cong
U(\g)\ten_{U(F^k\n)}{\mathcal S}_{F^{k+1}\n}^{F^k\n}\ten_{U(\n^k)}U(\n^k)^*
\cong
U(\g)\ten_{U(F^k\n)}U(F^k\n)^*={\mathcal S}_{F^k\n}.
$$
The Lemma is proved.
\end{proof}

The statement of the Theorem follows immediately from the previous Lemma.
\end{proof}

Let us return to the setting from the beginning of the section so that
$\g=\on{Lie}(G)$ for a semisimple group $G$ etc. 
By the previous construction for any $w\in W$ we obtain the semiregular $\g$-bimodule
${\mathcal S}_{\n_w^-}$ integrable over ${\n_w^-}\oplus{\n_w^-}$

\subsection{Semiregular bimodules over a semisimple Lie algebra. Local cohomology realization.}
Below we present a geometric description of the modules ${\mathcal S}_{\n_w^-}$ in terms of 
local cohomology of the structure sheaf $\CO_G$ on $G$ with support in the nilpotent subgroups 
$N_w$ such that $n_w^-=\on{Lie}N_w$. 

Let us first recall basic facts concerning local cohomology of quasicoherent sheaves.

\subsubsection{Local cohomology of quasicoherent sheaves.} Below we use the definition of local 
cohomology that differs from the original one but is equivalent to it. 

\begin{defn}
Consider a quasicoherent sheaf $\mathcal F$ on a scheme $X$ with a fixed locally closed subscheme
$V\subset X$. Denote the formal scheme obtained by completion of $V$ in 
$X$ by $V_X$. In particular we have an embedding $j:\  V_X\hookrightarrow X$. 
We define local cohomology spaces of 
$\mathcal F$ with support in $V$ as follows:
$$
H^\bullet_V(X,{\mathcal F}):=\on{R}\Gamma(V_X,j^!({\mathcal F})).
$$
Here the functor $j^!( {\mathcal F})$ is by definition the derived functor of the functor 
${\mathcal H}om_{{\mathcal O}_X}({\mathcal O}_{V_X},{\mathcal F})$.
\end{defn}

In particular denote the projection $X\to pt$ by $\pi$. Then we have
$H^\bullet_V(X,{\mathcal F})=H^\bullet(\on{R}\pi_*j^!({\mathcal F}))$.
\vskip 1mm
\Rem
Note that the definition is independent of the characteristic of the basic field and can be used 
in the prime characteristic as well. 

\subsubsection{Local cohomology of the structure sheaf} Here we collect some
 nessesary facts about local cohomology of the structure sheaf to be used later.

\begin{lem} \label{propert}

\vskip 0.25 mm\noindent\begin{itemize}
\item[(i)] In characteristic zero the  algebra of global differential operators $D_X$ acts on the spaces 
$H^\bullet_V(X,{\mathcal O}_X)$. 

\item[(ii)]
Let $V\subset X$ be a smooth closed subvariety of a smooth variety of codimension $d$.
The characteristic of the base field is supposed to be $0$. Then $H^{\ne d}_V(X,{\mathcal O}_X)=0$
and  $H^{ d}_V(X,{\mathcal O}_X)$ is isomorphic to the global sections of the D-module 
$i^V_!({\mathcal O}_V)$. Here $i^V$ denotes the embedding of $V$ into $X$.
\qed
\end{itemize}
\end{lem}

Note in particular that given an action of an algebraic group $G$ on $X$ we obtain the map
of Lie algebras $\g:=\on{Lie}(G)\map {\mathcal V}ect(X)$. Thus (in characteristic zero case)
the Lie algebra $\g$ acts on the local cohomology spaces $H^\bullet_V(X,{\mathcal O}_X)$.

\subsubsection{Local cohomology realization of semiregular bimodules.}
Let us return to the setup of the beginning of the section and consider the $\g$-$\g$-bimodule 
${\mathcal S}_{\n_w^-}$ for an element of the Weyl group $w\in W$.
Using the above Lemma we obtain a nice geometric description for ${\mathcal S}_{\n_w^-}$
 in terms of local cohomology of the structure sheaf on the group $G$.

Consider the space of local cohomology of ${\mathcal O}_G$ with support in the nilpotent subgroup
$N_w\subset G$. Recall that $\on{Lie}(N_w)=\n^-_w=\n^-\cap w(\n^+)$.
Note that $G$ acts on itself both by left and right translations and this provides two inclusions
$\g\map{\mathcal V}ect(G)$.

\begin{lem}
The $\g$-$\g$-bimodules $H^{ \dim G-\ell(w)}_{N_w}(G,{\mathcal O}_G)$ and 
 ${\mathcal S}_{\n_w^-}$ are naturally isomorphic. The spaces $H^{ \ne \dim G-\ell(w)}_{N_w}
(G,{\mathcal O}_G)$ vanish.
\end{lem}

\begin{proof}
It follows from the definition of the functor 
$i_!^{N_w}$ that $\Gamma(i_!^{N_w}({\mathcal O}_{N_w}))$
is isomorphic to $U(\g)\ten_{U(\n_w^-)}{\mathcal O}_{N_w}=U(\g)\ten_{U(\n_w^-)}\left(
{U(\n_w^-)}\right)^*={\mathcal S}_{\n_w^-}$.
Now the  statement of the Lemma follows form Lemma~\ref{propert}(ii).
\end{proof}
\vskip 1mm
\noindent
\Rem
In the next section we will prove that the statement of the  Lemma remains true in the 
positive characteristic case. Yet the proof in that case becomes more complicated and uses 
deep properties of local cohomology.

\subsection{Twisting functors and twisted Verma modules.}
Let us introduce the categories of $\g$-modules we will work with.

Denote by $\g\mod$ the category of $\h$-integrable $\g$-modules $$M=
\underset{\lambda \in X}\bigoplus M_\lambda, \dim M_\lambda <\infty, $$ where $M_\lambda$
 denotes the weight space of the weight $\lambda$, i.e. $M_\lambda=\{m\in M|\text{ for any }
h\in\h:\ h(m)=\lambda(h)m\}$. 

Let ${\mathcal O}(\g)$ be the full subcategory of $\g\mod$ consisting of $\n^+$-locally finite
$\g$-finitely generated modules.

Below for any $w\in W$ we define a functor $S_w:\ \g\mod\map\g\mod$ and a functor 
$\Theta_w:\ {\mathcal O}(\g)\map{\mathcal O}(\g)$. 

\subsubsection{Twisting functors on the category $\g\mod$.} The variant of a twisting 
functor easiest to define is as follows.

\begin{defn}
For $M\in\g\mod$ consider the $\g$-module $$S_w(M):={\mathcal S}_w\ten_{U(\g)}M.$$
\end{defn}

\begin{lem}
The functor $S_w$ is well defined on the category $\g\mod$. In other words for $M\in\g\mod$ 
the module  $S_w(M)$ is $\h$-diagonizible with finite dimensional $X$-grading components.
\qed
\end{lem}
 Our goal however is to define a similar functor on the category ${\mathcal O}(\g)$. The functor 
 $S_w$ evidently does not preserve the category, more precisely it is easy to verify that $S_w$ takes
$\n^+$-integrable modules to $w^{-1}(\n^+)$-integrable modules. Thus we need to add a twist of the
$X$-grading.

Recall that the group $W$ equals the quotient of the normalizer of $T$ in $G$ by $T$ itself. For any
$w\in W$ we choose its representative ${\mathfrak w}\in\on{Norm}(T)\subset G$. 

\begin{defn}
We define the functor of the grading twist $T_w:\ \g\mod\map\g\mod$ by setting
$T_w(M)_\lambda:=M_{w(\lambda)}$ with the $\g$-action given by
$$
g\cdot m:= Ad_{\mathfrak w}(g)( m)\text{ for } m\in M, g\in \g.
$$
The functor $\Theta_w:\ \g\mod\map\g\mod$ is defined as the composition
$\Theta_w=T_w\circ S_w$.
\end{defn}

\begin{lem}
\vskip 0.25 mm\noindent\begin{itemize}
\item[(i)]
The functor $\Theta_w$ is well defined as a functor ${\mathcal O}(\g)\map{\mathcal O}(\g)$.
\item[(ii)]
$\Theta_w:\ {\mathcal O}(\g)\map{\mathcal O}(\g)$ is right exact.
\end{itemize}
\end{lem}

We extend the functor $\Theta_w$ to the derived category ${\mathsf D}^b
({\mathcal O}(\g))$ in the usual way.
Abusing notation we denote the obtained functor by the same letter.

\begin{prop}
The functor  $\Theta_w:\ {\mathsf D}^b({\mathcal O}(\g))\map {\mathsf D}^b
({\mathcal O}(\g))$ is an equivalence of
 the triangulated categories.
\end{prop}

\begin{proof}
Essentially one has to construct the quasi-inverse functor for 
${\mathcal S}_w\overset{\on{L}}\ten\cdot$ on suitably chosen derived categories.
Consider the functor $\g\mod\map\g\mod$ as follows:
$$
M\mapsto \left(\hom_\g({\mathcal S}_w,M)\right)^{T\on{-ss}},
$$
where $(\cdot)^{T\on{-ss}}$ denotes taking the maximal $T$-semisimple submodule.

\begin{lem}
\vskip 0.25mm
\noindent
\begin{itemize}
\item[(i)] The functor $\left(\hom_\g({\mathcal S}_w,\cdot)\right)^{T\on{-ss}}$ is well defined on the 
category $\g\mod$ and is left-exact.
\item[(ii)] The functor $\left(\hom_\g({\mathcal S}_w,\cdot)\right)^{T\on{-ss}}\circ T_{w^{-1}}$
 is well defined as a functor ${\mathcal O}(\g)\map{\mathcal O}(\g)$.\qed
\end{itemize}
\end{lem}

Denote the right derived functor of $\left(\hom_\g({\mathcal S}_w,\cdot)\right)^{T\on{-ss}}
\circ T_{w^{-1}}$ by
$$
\underline{\on{R}\hom}\bul_\g({\mathcal S}_w,\cdot):\ {\mathsf D}^b
({\mathcal O}(\g))\map
{\mathsf D}^b
({\mathcal O}(\g)).
$$
Consider also the functor on the category of $\g$-bimodules $M\mapsto
\left(\hom_\g({\mathcal S}_w,M)\right)^{\on{ad}_\g}$, where $(\cdot)^{\on{ad}_\g}$
denotes the maximal $\on{ad}_\g$-locally finite submodule in $(\cdot)$. 

\begin{lem}\vskip 0.25mm
\noindent
\begin{itemize}
\item[(i)]
The functor $\left(\hom_\g({\mathcal S}_w,\cdot)\right)^{\on{ad}_\g}$ is left-exact.

\item[(ii)] We have $\left(\hom_\g({\mathcal S}_w,{\mathcal S}_w)\right)^{\on{ad}_\g}\til\map
U(\g)$.
\qed
\end{itemize}
\end{lem}

For any $M\bul\in {\mathsf D}^b({\mathcal O}(\g))$ we have a natural map
\begin{gather*}
\tau:\ M\til\map\hom_{\g}({\mathcal S}_w,{\mathcal S}_w)^{\on{ad}_\g}
\overset{\on{L}}\ten M\\ \map
\on{R}\bul\left(\hom_{\g}({\mathcal S}_w,{\mathcal S}_w)^{\on{ad}_\g}\right)
\overset{\on{L}}\ten M\map
\underline{\on{R}\hom}\bul_{\g}({\mathcal S}_w,{\mathcal S}_w\overset{\on{L}}\ten M).
\end{gather*}
\begin{lem}
The map $\tau$ is an isomorphism in the category ${\mathsf D}^b({\mathcal O}(\g))$.
\end{lem}
\begin{proof}
First we apply the functor 
$\underline{\on{R}\hom}\bul_{\g}({\mathcal S}_w,{\mathcal S}_w\overset{\on{L}}\ten \cdot)$
to a Verma module $M(\lambda)$.
On the level of vector spaces (not regarding the $\g$-module structure) we have
\begin{gather*}
\underline{\on{R}\hom}\bul_{\g}({\mathcal S}_w,{\mathcal S}_w\overset{\on{L}}\ten M(\lambda))
\til\map
\underline{\on{R}\hom}\bul_{\n_w}((U(\n_w))^*,(U(\n_w))^*\ten U(\n^-\cap w(\n^-)))\\
\til\map\left(\Hom_{\mathbb C}((U(\n_w))^*,U(\n^- \cap w(\n^-)))\right)^{T\on{-ss}}\til\map U(\n^-).
\end{gather*}
We leave to the reader to check that the above module is isomorphic to $M(\lambda)$ and that $\tau$
provides this isomorphism. 

It follows that for any $\g$-module $M$ from the category ${\mathcal O}(\g)$ the map $\tau$ provides an isomorphism
$M\til\map\underline{\on{R}\hom}\bul_{\g}({\mathcal S}_w,{\mathcal S}_w\overset{\on{L}}\ten M)$.

Now note that for any finite complex of modules from the category ${\mathcal O}(\g)$ is quasi-isomorphic to 
a finite complex of modules finitely filtered by Verma modules. This assertion implies the 
statement of the Lemma.
\end{proof}
Now the statement of the proposition follows from the Lemma since the functor
$T_{w^{-1}}$ is quasi-inverse for the functor $T_w$.
\end{proof}

\subsubsection{Twisted Verma modules}
As usual, we define the Verma module we with the highest weight $\lambda$ by
$M(\lambda):=U(\g)\ten_{\h\oplus\n^+}\BC(\lambda)$, where $\BC(\lambda)$ denotes the one 
dimensional weight module over $\h$.

\begin{defn}
We define the twisted Verma module with the twist $w\in W$ and the highest weight $\lambda$\
by $M^w(\lambda)=\Theta_w(M(w^{-1}\cdot\lambda))$. Here as usual the dot action of the 
Weyl group on the weights is given by $w\cdot \lambda:=w(\lambda+\rho)-\rho$.
\end{defn}

Recall that the character of a module $M$ from the category ${\mathcal O}(\g)$ is defined as the
formal expression $\ch(M)=\sum(\dim M_\lambda) e^\lambda$. In particular the character of the 
Verma module $M(\lambda)$
$$
\ch(M(\lambda))=e^\lambda\frac1{\prod_{\alpha\in R^+}(1-e^{-\alpha})}.
$$
\begin{lem}
The character of a twisted Verma module $M_w(\lambda)$ coincides with the one of the 
Verma module $M(\lambda)$.
\qed
\end{lem}
\vskip 1mm
\noindent
\Rem
Still the twisted Verma module $M_w(\lambda)$ usually differs from the Verma module with the 
same highest weight. For example for $G=SL(2)$ the twisted Verma module with the only 
possible nontrivial twist and the highest weight $0$ is isomorphic to the {\em contragradient}
Verma module of the highest weight $0$.
  
\subsection{Contragradient quasi-Verma modules in characteristic $0$} 
Recall that the category ${\mathcal O}(\g)$
is decomposed into the direct sum of linkage classes or blocks. We will be interested in the rough 
decomposition $ {\mathcal O}(\g)=\underset{\theta\in\on{Spec}\CZ} {\mathcal O}(\g)_\theta$. Here
$\CZ$ denotes the center of the universal enveloping algebra of $\g$ and the subcategory
${\mathcal O}(\g)_\theta$ consists of the modules $M\in  {\mathcal O}(\g)$ such that for any 
$m\in M$ and $z\in\CZ$ there exists an integer $k$ such that $(z-\theta(z))^k(m)=0$. In particular
consider the block  $ {\mathcal O}(\g)_0$ corresponding to the trivial central character. 

It is known that $M(\lambda)$ belongs to $ {\mathcal O}(\g)_0$ if and only if the weight $\lambda$ 
is of the form $\lambda=w\cdot 0$ for some $w\in W$. Let $w_0$ be the longest element in the 
Weyl group. It is known also that $M(w_0\cdot 0)$ is simple. 

More generally, consider a regular dominant integral weight $\lambda$ and the block 
${\mathcal O}(\g)_\theta$ containing $M(\lambda)$. Then it is known that $M(\mu)$ 
belongs to $ {\mathcal O}(\g)_\theta$ if and only if the weight $\mu$ 
is of the form $\mu=w\cdot \lambda$ for some $w\in W$.
 Again it  is known  that $M(w_0\cdot \lambda)$ is simple.

\begin{defn}
Fix the block ${\mathcal O}(\g)_\theta$ containing $M(\lambda)$ with regular integral dominant 
highest weight $\lambda$. Let $\mu=ww_0\cdot\lambda$. Then the contragradient 
quasi-Verma module with the 
highest weight $\mu$ is defined as follows:
$$\D\til{M}(\mu):=\Theta_w(M(w_0\cdot\lambda)).$$
\end{defn}
\vskip 1mm
\noindent
\Rem
By the very definition a contragradient quasi-Verma module provides a  special  
case of twisted Verma module with 
a certain choice of the highest weight and the twisting element of the Weyl group. Moreover, the
 next statement shows that in characteristic $0$ the class of contragradient 
quasi-Verma modules has a complete 
and nice description. 

\begin{prop}
For any regular dominant integral weight $\lambda$ and any $w\in W$ the 
contragradient quasi-Verma module 
$\D\til{M}(\mu)$ with the highest weight $\mu=ww_0\cdot\lambda$ is isomorphic to the contragradient 
Verma module $\D M(\mu)$.
\end{prop}

\begin{proof}
First note that there exists a canonical map $\pi:\ \D\til{M}(\mu)\map\D M(\mu)$ since 
$\D\til{M}(\mu)$ 
has a highest weight vector of the weight $\mu$. Now apply the functor $(\Theta_w)^{-1}$ in 
the derived category to the modules $\D\til{M}(\mu)$, $\D M(\mu)$ and to the map $\pi$.
Evidently by the definition of a contragradient quasi-Verma module we have 
$(\Theta_w)^{-1}\D\til{M}(\mu)=M(w_0\cdot\lambda)$.

\begin{lem}
The complex $(\Theta_w)^{-1}(\D M(\mu))$ is concentrated in homological degree $0$. Moreover
the character of  $(\Theta_w)^{-1}(\D M(\mu))$ coincides with the one of $M(w_0\cdot\lambda)$.
\qed
\end{lem}

Since the module  $M(w_0\cdot\lambda)$ is simple, we conclude that the map $(\Theta_w)^{-1}(\pi)$
is an isomorphism in the derived category. Thus it is an isomorphism in ${\mathcal O}(\g)$. It 
follows that $\pi$ is an isomorphism itself.
\end{proof}
\subsubsection{Bernstein-Gelfand-Gelfand complex.} We conclude this section with the well-known 
result. In fact in the next two sections we will be dealing with direct characteristic $p$ analogues of it.
\begin{prop}
Fix a regular dominant integral weight $\lambda$. Then there exists a complex of $\g$-modules
$K_{\BC}\bul(\lambda)$ (called the contragradient Bernstein-Gelfand-Gelfand complex) as
follows:
$$ K_{\BC}^m(\lambda)=\underset{\ell(w)=m}\bigoplus\ \D M(w\cdot \lambda).
$$
Higher cohomology spaces of the complex vanish and $H^0(K_{\BC}\bul(\lambda))$ is isomorphic to
the simple $\g$-module with the highest weight $\lambda$. \qed
\end{prop}

\section{Positive characteristic case} 
\subsection{Semiregular modules over a semisimple Lie algebra. Algebraic setting in positive 
characteristic.}
>From now on we work over the algebraic closure $\F$ of the finite field $\f$ of characteristic $p$.
Again we fix the root data $(Y,X,\ldots)$ of the finite type $(I,\cdot)$. Consider the algebraic group 
$G_{\F}$ over the field $\F$ corresponding to the root data.

Denote by $\U_{\F}(\g)$ the specialization of the Kostant integral form of the universal enveloping 
algebra for $\g$ specialized at $\F$.

Our aim now is to define semiregular modules in the new setting. Consider the subalgebra
$\U_{\F}(\n_w^-)\subset \U_{\F}(\g)$ and its natural bimodule ${\mathcal O}_{N_{w,\F}}$.

\begin{defn}
The  $\U_{\F}(\g)$-module ${\mathcal S}_{\n_w^-,\F}:=\U_{\F}(\g)\ten_{\U_{\F}(\n_w^-)}
{\mathcal O}_{N_{w,\F}}$ is called the left semiregular module corresponding to the element
$w$ of the Weyl group.
\end{defn}
\vskip 1mm
\noindent
\Rem
Evidently the definition simply mimics the one in characteristic zero and it is not clear a priori that
the defined object makes sense. In particular it is much more difficult in the present setup to define the 
$\U_{\F}(\g)$-bimodule structure on  ${\mathcal S}_{\n_w^-,\F}$.

\subsection{More local cohomology}
To move further we will need a few general results on local cohomology.

Let us mention first that for any variety $X$ over $\F$ equipped with an action of the group
$G_{\F}$ and for any locally closed subscheme $V\subset X$
the hyperalgebra $\U_{\F}(\g)$ acts on the local cohomology spaces 
$H^{\bullet}_V(X,{\mathcal O}_X)$. This follows from Lemma~\ref{propert}(ii).

Suppose we have a smooth surjective morphism of smooth schemes $p:\ X\map Y$ of relative 
dimension  $d$. Suppose that the relative canonical bundle is trivial. Let $V$ be a closed subscheme of 
$X$ projecting isomorphically to $W\subset Y$. Denote the embedding of the formal completion of
$V$ in $X$ (resp. of $W$ in $Y$) by $\til{j}:\ V_X\map X$ (resp. by $j:\ W_Y\map Y$).

\begin{prop} \label{smooth}
There exist the canonical maps of local cohomology spaces
$$
H^{\bullet}_V(X,{\mathcal O}_X)\map
H^{\bullet-d}_W(Y,{\mathcal O}_Y)
$$
\end{prop}

\begin{proof} Denote by $\pi$ the projection $W_Y\map pt$. Then by definition of local cohomology 
we have 
\begin{gather*}
H^{\bullet}_V(X,{\mathcal O}_X)=H\bul(\on{R}\pi_*\on{R}p_*\til{j}^!{\mathcal O}_X)\\
\til\map H\bul(\on{R}\pi_*\on{R}p_*\til{j}^!p^*{\mathcal O}_Y)\til\map
H\bul(\on{R}\pi_*\on{R}p_*\til{j}^!p^!{\mathcal O}_Y[d])\\ \til\map
H\bul(\on{R}\pi_*\on{R}p_*p^!j^!{\mathcal O}_Y[d]).
\end{gather*}
Here we used that the map $p$ is smooth of relative dimension $d$. Next note that the map
$p:\ V_X\map W_Y$ is pro-finite, thus we can replace $p_*$ by $p_!$ in the above formula.
We obtain the following map
\begin{gather*}
H^{\bullet}_V(X,{\mathcal O}_X)
\til\map
H\bul(\on{R}\pi_*\on{R}p_!p^!j^!{\mathcal O}_Y[d])\map H\bul(\on{R}\pi_*
j^!{\mathcal O}_Y[d]) =H^{\bullet-d}_W(Y,{\mathcal O}_Y). 
\end{gather*}
Here we used the canonical adjointness  map $p_!p^!\map\on{Id}$.
\end{proof}

\subsubsection{Definition of the hyperalgebra for $\g$ via local cohomology.}
Consider the action of the group $G_{\F}$ on a variety $X$ over $\F$ such that the action map
$X\times G_{\F}\map X$ is smooth. Let $V\subset X$ be a 
locally closed subscheme. Denote by $\{e\}\subset G_{\F}$ the unit element.
\begin{lem}
\vskip 0.25 mm\noindent\begin{itemize}
\item[(i)]
There is a natural map of local cohomology spaces
$$
H\bul_{V\times\{e\}}(X\times G_{\F},{\mathcal O}_{X\times G_{\F}})\map
H^{\bullet-\dim G}_{V}(X,{\mathcal O}_X).
$$
\item[(ii)]
In particular the multiplication map $G_{\F}\times G_{\F}\map G_{\F}$ makes
$$H^{\dim G}_{\{e\}}(G_{\F},{\mathcal O}_{G_{\F}})$$
into an associative algebra.
\item[(iii)]
The canonical map from the first part of the Lemma makes 
$H^{\bullet}_{V}(X,{\mathcal O}_X)$ a $H^{\dim G}_{\{e\}}(G_{\F},
{\mathcal O}_{G_{\F}})$-module.
\end{itemize}
\end{lem}
\begin{proof}
All the assertions of the Lemma follow immediately from Proposition~\ref{smooth}.
\end{proof}
\begin{defn}
The algebra $H^{\dim G}_{\{e\}}(G_{\F},{\mathcal O}_{G_{\F}})$ from the previous Lemma is 
called the hyperalgebra for $\g$ over $\F$.
\end{defn}
The following statement is well known.
\begin{prop} \label{2alg}
The algebra $H^{\dim G}_{\{e\}}(G_{\F},{\mathcal O}_{G_{\F}})$ is naturally isomorphic to the
specialization $\U_{\F}(\g)$ of the Kostant integral form for the universal enveloping algebra for $\g$.
\qed
\end{prop}
\vskip 1mm
\noindent
\Rem
In particular suppose there is a smooth action of the group $G_{\F}$ on a variety $X$.
Then the algebra $\U_{\F}(\g)$ acts on any local cohomology spaces $H\bul_{V}(X,{\mathcal O}_X)$
(for any locally closed subscheme $V\subset X$).

\begin{cor} \label{action}
(from Proposition~\ref{smooth}  and Proposition~\ref{2alg})
Suppose that an affine algebraic group $G_{\F}$ of dimension $d$ acts freely on a smooth variety $X$ 
defined over $\F$ with the quotient $Y=X/G_{\F}$. 
Suppose also that a closed subvariety $V\subset X$ maps isomorphically onto its image under the 
projection $W\subset Y$.
Then there exists a canonical map
$$
\left(H^{\bullet}_V(X,{\mathcal O}_X) \right)_{\U_{\F}(\g)\on{-coinv}}\map
H^{\bullet-d}_W(Y,{\mathcal O}_Y).
\qed
$$
\end{cor}

Let $V$ be a smooth subvariety in a smooth variety $X$. As usual denote by $T^*_VX$ the conormal
 bundle to $V$. We have an embedding $V\hookrightarrow T^*_VX$ via the zero section.

\begin{prop} \label{filt}
There exists a natural filtration $F\bul$ on local cohomology spaces  
$H^{\bullet}_V(X,{\mathcal O}_X)$ with the associated graded spaces 
$\gr^{F\bul}\left(H^{\bullet}_V(X,{\mathcal O}_X)\right)$ isomorphic to 
$H^{\bullet}_V(T^*_VX,{\mathcal O}_{T^*_VX})$.
\end{prop}

\begin{proof}
Denote by $V^{(n)}_X$ the $n$-th infinitesimal neighborhood of $V$ in $X$. Let
$j^{(n)}:\ V^{(n)}_X\hookrightarrow X$. Then we have a natural filtration $F\bul$ of
$H\bul(\on{R}\Gamma(V_X,j^!{\mathcal O}_X))$ by 
$H\bul(\on{R}\Gamma(V^{(n)}_X,j^{(n)!}{\mathcal O}_X))$.

\begin{lem}
The spaces $\gr^{F\bul}\left(H^{\bullet}_V(X,{\mathcal O}_X)\right)$ are isomorphic to 
$H^{\bullet}_V(T^*_VX,{\mathcal O}_{T^*_VX})$ with the second grading given by the 
action of the multiplicative group of $\F$.
\qed
\end{lem}

\end{proof}

\subsection{Semiregular modules over the hyperalgebra for a semisimple Lie group
 in characteristic $p$. Local cohomology realization.}
Here at last we introduce the $\U_{\F}(\g)$-bimodule structure on ${\mathcal S}_{n_w^-,\F}$.

Note that by Lemma~\ref{propert}(ii) the algebra of differential operators with divided powers 
on $G_{\F}$  acts on $H\bul_{N_{w,\F}}(G_{\F},{\mathcal O}_{G_{\F}})$. On the other hand
the group $G_\F$ acts on itself both by left and right translations and this provides two embeddings
of $\U_{\F}(\g)$ into  the algebra of differential operators with divided powers on $G_{\F}$. Like in 
the complex case our aim is to identify the $\U_{\F}(\g)$-bimodules ${\mathcal S}_{n_w^-,\F}$
and $H\bul_{N_{w,\F}}(G_{\F},{\mathcal O}_{G_{\F}})$.

Consider the map $\mu:\ G_{\F}\times N_{w,\F}\map G_{\F}$ provided by the multiplication. One 
can view this map as taking quotient of the variety $G_{\F}\times N_{w,\F}$ by the free diagonal 
action of the group $N_{w,\F}$.

\begin{lem}
There is a canonical map of the left $\U_{\F}(\g)$-modules
$$m: \U_{\F}(\g)\ten_{\U_{\F}(\n_w^-)}{\mathcal O}_{N_{w,\F}}
\map H^{\dim G-\dim N_w}_{N_{w,\F}}(G_{\F},{\mathcal O}_{G_{\F}}).$$
\end{lem}

\begin{proof}
First note that 
\begin{gather*}
\U_{\F}(\g)\ten_{\U_{\F}(\n_w^-)}{\mathcal O}_{N_{w,\F}}\til\map
\left(H^{\dim G}_{\{e\}}(G_{\F},{\mathcal O}_{G_{\F}})\ten H_{N_{w,\F}}^0(N_{w,\F},
{\mathcal O}_{N_{w,\F}})\right)_{\U_{\F}(\n_w^-)\on{-coinv}}\\ \til\map
\left(H^{\dim G}_{\{e\}\times N_{w,\F}}(G_{\F}\times N_{w,\F},
{\mathcal O}_{G_{\F}\times N_{w,\F}})\right)_{\U_{\F}(\n_w^-)\on{-coinv}}.
\end{gather*}
Now consider the map $\mu:\ G_{\F}\times N_{w,\F}\map G_{\F}$ and use Corollary~\ref{action}.
The obtained morphism is the one of $\U_{\F}(\g)$-modules because the map $\mu$ is equivariant
with respect to the action of $G_{\F}$ by left translations.

\end{proof}

\begin{prop}
The canonical map from the previous Lemma is an isomorphism of the left $\U_{\F}(\g)$-modules.
\end{prop}

\begin{proof} 
Let us identify the normal bundle to $\{e\}\times N_{w,\F}$ in $G_{\F}\times N_{w,\F}$
(resp. the  normal bundle to $ N_{w,\F}$ in $G_{\F}$) with $\g_{\F}\times  \n^-_{w,\F}$
(resp. with $\g_{\F}$). By the previous Lemma we have the map of the filtered objects
$$
\left(H^{\dim G}_{\{e\}\times N_{w,\F}}(G_{\F}\times N_{w,\F},
{\mathcal O}_{G_{\F}\times N_{w,\F}})\right)_{\U_{\F}(\n_w^-)\on{-coinv}}\map
H^{\dim G-\dim N_w}_{N_{w,\F}}(G_{\F},{\mathcal O}_{G_{\F}}).
$$
Now use Proposition~\ref{filt} and consider the morphism of associated graded objects
\begin{multline*}
\gr^{F\bul} m:\ \left(H^{\dim \g}_{\{0\}}(\g_{\F},{\mathcal O}_{\g_{\F}})\ten
H^0(N_{w,\F},{\mathcal O}_{N_{w,\F}}\right)_{\U_{\F}(\n_w^-)\on{-coinv}}\\
\map H^{\dim \g-\dim\n_w^-}_{\n_{w,\F}^-}
(\g_{\F},{\mathcal O}_{\g_{\F}}).
\end{multline*}
\begin{lem}
The morphism $\gr^{F\bul} m$ is an isomorphism of vector spaces.
\qed
\end{lem}
\end{proof}

We have proved the following statement.

\begin{thm} \label{local}
For any $w\in W$ the  semiregular $\U_{\F}(\g)$-module is isomorphic to 
$H^{\dim G-\dim N_w}_{N_{w,\F}}(G_{\F},{\mathcal O}_{G_{\F}})$ and is equipped with the 
natural $\U_{\F}(\g)$-bimodule structure.
\qed
\end{thm}

Just like in the characteristic $0$ case, we use the semiregular bimodules to construct certain functors
on the category of $\U_{\F}(\g)$-modules.

\subsection{Twisting functors and twisted Verma modules in characteristic $p$.}
Below we work in the category of $\U_{\F}(\g)$-modules integrable over the Cartan subgroup 
$T_{\F}$.
 We denote this category by $\U_{\F}(\g)\mod$. 

\begin{defn}
For $M\in\U_{\F}(\g)\mod$ consider the $\U_{\F}(\g)$-module 
$$S_{w,\F}(M):={\mathcal S}_{w,\F}\ten_{\U_{\F}(\g)}M.$$
\end{defn}

\begin{lem} The functor $S_{w,\F}$ is well defined on the category $\U_{\F}(\g)\mod$. In other 
words it takes $T_{\F}$-integrable $\U_{\F}(\g)$-modules to $T_{\F}$-integrable 
$\U_{\F}(\g)$-modules.
\qed
\end{lem}

Again we have to add a certain grading twist to this functor. Namely, as before, for $w\in W$
consider a lift of the element of the Weyl group to an element ${\mathfrak w}\in G_{\F}$.

\begin{defn}
We define the functor of the grading twist $T_w:\ \U_{\F}(\g)\mod\map\U_{\F}(\g)\mod$:
$T_w(m)=M$ as a vector space but
$$
u\cdot m:=\on{Ad}_{  \mathfrak w}(u)\cdot m\text{ for } u\in \U_{\F}(\g),\ m\in M.
$$
The functor $\Theta_{w,\F}$ is defined as the composition: $\Theta_{w,\F}:=T_{w,\F}\circ S_{w,\F}$.
\end{defn}
\vskip 1mm
\noindent
\Rem
We do not know any suitable analogue of the category $\mathcal O$ in characteristic $p$ such that
the functor $\Theta_{w,\F}$ would extend to its derived category and provide an autoequivalence 
of the triangulated category. We doubt the very possibility of a statement of this nature since the picture 
in characteristic $p$ becomes more complicated. This will be discussed below.

Consider the Verma module $M_{\F}(\lambda):=\U_{\F}(\g)\ten_{\U_{\F}(\h\oplus\n^+)}
\F(\lambda)$ with the highest weight $\lambda$.

\begin{lem}
Like in the characteristic zero case, the character of the module $\Theta_{w,\F}(M_{\F}(
w^{-1}\cdot\lambda))$ coincides with the one of $M_{\F}(\lambda)$.
\qed
\end{lem}

\begin{defn}
We call the module  $\Theta_{w,\F}(M_{\F}(
w^{-1}\cdot\lambda))$ the twisted Verma module with the twist $w\in W$ and the highest weight
$\lambda$ and denote it by $M_{\F}^w(\lambda)$. Let $\lambda$ be a regular dominant integral 
weight.
We call the module  $\Theta_{w,\F}(M_{\F}(w_0\cdot\lambda))$ the contragradient quasi-Verma module
with the highest weight $\mu=ww_0\cdot\lambda$ and denote it by $\D\til{M}_{\F}(\mu)$.
\end{defn}

\subsection{quasi-Verma modules in characteristic $p$ and the Grothendieck-Cousin complex
of $(G/B)_{\F}$} Below we explain how contragradient quasi-Verma modules appear naturally in the geometry of 
the 
Flag variety in characteristic $p$. 

Consider the Flag variety $\B_{\F}:=G_{\F}/B_{\F}.$ Here $B_{\F}=B^+_{\F}$ 
denotes the standard positive Borel subgroup in $G_{\F}$. Let us recall the facts concerning 
Bruhat decomposition of $\B_{\F}$.

\begin{lem}
\vskip 0.25 mm\noindent\begin{itemize}
\item[(i)]
The orbits of $B_{\F}$ on $\B_{\F}$ are enumerated by the elements of the Weyl group:
$$
\B_{\F}=\underset{w\in W}{\bigsqcup}C_{w,\F},
$$
where the orbit $C_{w,\F}=B_{\F}{\mathfrak w}B_{\F}/B_{\F}$.
\item[(ii)]
The orbit $C_{w,\F}$ is isomorphic to the $\F$-affine space of the dimension equal to 
$\dim \B_{\F}-\ell(w)=\sharp(R^+)-\ell(w)$.\qed
\end{itemize}
\end{lem}

The  $B_{\F}$-orbits on the Flag variety are called the Schubert cells.

\subsubsection{Global Grothendieck-Cousin complexes on $\B_{\F}$.}
Recall that the contragradient Weyl module with the regular integral dominant weight $\lambda$
is defined by
$$
\D W_{\F}(\lambda)=\left(\hom_{\U(\h\oplus\n^-)}(\U(\g),\F(\lambda))\right)^{\on{int}},
$$
where $(*)^{\on{int}}$ denotes the maximal $G_{\F}$-integrable submodule.

Consider the standard equivariant liner bundle $\CL(\lambda)$ on $\B_{\F}$. The following statement 
is known as the Borel-Weyl-Bott theorem in positive characteristic and is due to Kempf.

\begin{prop}
Higher cohomology groups of $\CL(\lambda)$ vanish and $$H^0(\B_{\F},\CL(\lambda))=
\D W_{\F}(\lambda).\qed$$
\end{prop}

Recall the following statement due to Kempf. 

\begin{thm}
\vskip 0.25 mm\noindent\begin{itemize}
\item[(i)] There exists a complex of $\U_{\F}(\g)$-modules $K_{\F}\bul(\lambda)$ with
$$
K_{\F}^m(\lambda)=\underset{\ell(w)=m}\bigoplus H^{m}_{C_{ww_0,\F}}
(\B_{\F},\CL(\lambda)).
$$
\item[(ii)] Higher cohomology spaces of the complex vanish and 
$H^0(K_{\F}\bul(\lambda))=\D W_{\F}(\lambda)$.
\qed
\end{itemize}
\end{thm}

Our aim is to give the algebraic interpretation of the local cohomology spaces
$H^{\ell(w)}_{C_{w,\F}}(\B_{\F},\CL(\lambda))$ in terms of contragradient quasi-Verma modules 
in characteristic $p$.

\subsubsection{Local cohomology construction of contragradient quasi-Verma modules.} Recall that 
in characteristic $0$ case all contragradient quasi-Verma modules appeared to be quasiisomorphic to ordinary 
contragradient Verma modules with appropriate highest weights. Below we present an analog
of this statement in characteristic $p$.

\begin{thm} \label{compare}
Fix the regular dominant integral weight $\lambda$.
Then the $\U_{\F}(\g)$-modules $\D\til{M}(w\cdot \lambda)$ and 
$H^{\ell(w)}_{C_{ww_0,\F}}(\B_{\F},\CL(\lambda))$ are naturally isomorphic.
\end{thm}

\begin{proof}
For the sake of simplicity we prove the Theorem only for the weight $\lambda=0$ 
(and thus $\CL(\lambda)={\mathcal O}_{\B_{\F}}$).  The proof in the general case does not differ
much from the one presented below. We start with reformulating the statement of the Theorem.

Consider the action of the group $N_{w,\F}$ on $\B_{\F}$ and  the locally closed subset 
$N_{w,\F}\cdot\{e\}$. Here $\{e\}\subset\B_{\F}$ denotes the unique zero dimensional Schubert cell.

\begin{lem}
The $\U_{\F}(\g)$-modules $$H^{\ell(w)}_{C_{ww_0,\F}}(\B_{\F},{\mathcal O}_{\B_{\F}})
\text{ and }
T_{ww_0,\F}\left(H^{\ell(w)}_{N_{w,\F}\cdot\{e\}}(\B_{\F},{\mathcal O}_{\B_{\F}})\right)$$
are naturally isomorphic.
\end{lem}

\begin{proof}
To prove the Lemma note only that the natural action of the element $ww_0$ of the
 Weyl group on $\B_{\F}$ takes  $N_{w,\F}\cdot\{e\}$ to the Schubert cell $C_{w,\F}$.
\end{proof}

Now consider the map $p:\ G_{\F}\times \B_{\F}\map \B_{\F}$ provided by the action of $G_{\F}$
on the Flag variety. On the other hand $p$ can be considered as taking quotient of
$G_{\F}\times \B_{\F}$ by the free diagonal action of $G_{\F}$.

\begin{lem} \label{map}
There exists a natural map of the left $\U_{\F}(\g)$-modules
$$
m:\ \left(H^{\codim N_w}_{N_{w,\F}\times\{e\}}(G_{\F}\times \B_{\F},
{\mathcal O}_{G_{\F}\times \B_{\F}})\right)_{\U_{\F}(\g)\on{-coinv}}\map
H^{\codim N_{w,\F}\cdot \{e\}}_{N_{w,\F}\cdot\{e\}}(\B_{\F},{\mathcal O}_{\B_{\F}}).
$$
Here the coinvariants in the LHS are taken along the action of 
$\U_{\F}(\g)$ on $$H^{\codim N_w}_{N_{w,\F}\times\{e\}}(G_{\F}\times \B_{\F},
{\mathcal O}_{G_{\F}\times \B_{\F}})$$ provided by 
 the diagonal action of $G_{\F}$ on $G_{\F}\times \B_{\F}$.
\end{lem}

\begin{proof}
The  map in question is constructed using Corollary~\ref{action}. The obtained morphism 
respects the $\U_{\F}(\g)$-module structures on the LHS and RHS since the map $p$ is 
$G_{\F}$-equivariant.
\end{proof}
\vskip 1 mm
\noindent
\Rem
Note that we have 
$$H^{\codim N_w}_{N_{w,\F}\times\{e\}}(G_{\F}\times \B_{\F},
{\mathcal O}_{G_{\F}\times \B_{\F}})=H^{\codim N_w}_{N_{w,\F}}
(G_{\F},{\mathcal O}_{G_{\F}})\ten H^{\dim \B_{\F}}_{\{e\}}(\B_{\F},{\mathcal O}_{\B_{\F}})$$
by the Kunneth formula.

The following statement is well-known.

\begin{lem}
There exist  natural isomorphisms of the left $\U_{\F}(\g)$-modules
$$H^{\dim \B_{\F}}_{\{e\}}(\B_{\F},{\mathcal O}_{\B_{\F}})\til\map M_{\F}
(w_0\cdot 0)\text{  and }
H^0_{C_{w_0,\F}}(\B_{\F},{\mathcal O}_{\B_{\F}})\til\map\D M_{\F}(0).
\qed$$
\end{lem}

Thus Lemma~\ref{map} compared with Theorem~\ref{local}
 provides a natural $\U_{\F}(\g)$-module morphism
$$m: {\mathcal S}_{w,\F}\ten_{\U_{\F}(\g)}M_{\F}(w_0\cdot 0)\map 
H^{\ell(w)}_{N_{w,\F}\cdot\{e\}}
( \B_{\F},{\mathcal O}_{\B_{\F}}).$$

\begin{prop} \label{iso}
The map $m$ is an isomorphism of the $\U_{\F}(\g)$-modules.
\end{prop}

\begin{proof}
First it is known that $\ch  ( H^{\ell(w)}_{C_{w,\F}}(\B_{\F},{\mathcal O}_{\B_{\F}}))$ is equal to 
the character  of the Verma module $M_{\F}(w\cdot 0)$. Next it is easy to verify that the 
$\U_{\F}(\g)$-module
$T_{ww_0,\F}( {\mathcal S}_{w,\F}\ten_{\U_{\F}(\g)}M_{\F}(w_0\cdot 0))$ has the same character.

The rest of the proof is left to the reader.
\end{proof}

Proposition~\ref{iso} completes the proof of the Theorem. 

\end{proof}
\begin{cor}
For a regular dominant integral weight $\lambda$ the contragradient quasi-Verma
 module $\D\til{M}_{\F}(\lambda)$
is isomorphic to the contragradient Verma module $\D M_{\F}(\lambda)$.
\qed
\end{cor}

\section{Algebraic construction of the differentials in the Grothendieck-Cousin complex in
 characteristic $p$}
\ssn
Our aim is to construct the complex
$K_{\F}\bul(\lambda)$ with
$$
K_{\F}^m(\lambda)=\underset{\ell(w)=m}\bigoplus H^{m}_{C_{ww_0,\F}}
(\B_{\F},\CL(\lambda)).
$$
by purely algebraic methods. So far we have partly achieved this: we have proved that
$$
K_{\F}^m(\lambda)=\underset{\ell(w)=m}\bigoplus \D\til{M}_{\F}(w\cdot\lambda).
$$
Below we show how to express the components of the differentials 
$K_{\F}^m(\lambda)\map K_{\F}^{m+1}(\lambda)$ in purely algebraic terms. 

\subsubsection{Example: $SL(2)$.} The picture is particularly easy for $G=SL(2)$. Let $\lambda$
be a regular dominant integral weight.
By definition we have the canonical embedding 
$$i:\ \D W_{\F}(\lambda)\map \D M_{\F}(\lambda).$$

\begin{lem}
The cokernel of the map $i$ is isomorphic to the Verma module $M_{\F}(-\lambda-2\rho)$.
\end{lem}

\begin{proof}
This follows from the existence of the Grothendieck-Cousin complex for $\B_{\F}$ in this case
and from the identifications $H^0_{C_{s,\F}}(\B_{\F},\CL(\lambda))=
\D\til{M}(\lambda) =\D M_{\F}(\lambda)$ and
$H_{\{e\}}^1(\B_{\F},\CL(\lambda))=\D\til{M}_{\F}(s\cdot\lambda)=M_{\F}(-\lambda-2\rho)$.
\end{proof}
\vskip 1mm
\noindent
\Rem
No doubt the map $H^0_{C_{s,\F}}(\B_{\F},\CL(\lambda))\map H_{\{e\}}^1(\B_{\F},\CL(\lambda))$ 
is just the 
coboundary map in the Grothendieck-Cousin complex of local cohomology in the $SL(2)$ case. 
However we prefer to use both the  geometric  and the algebraic language here.
In particular the above Lemma can be thought of as an algebraic construction of the map
$$
d_{{\mathfrak {sl}}(2)}: \D\til{M}_{\F}^{{\mathfrak {sl}}(2)}(\lambda)\map
\D\til{M}_{\F}^{{\mathfrak {sl}}(2)}(-\lambda-2\rho).
$$
\sssn
Recall that the Weyl group $W$ is the Coxeter group with the set of generators $\{s_i|i\in I\}$.
Every element $w\in W$ has the shortest expression via the generators of the form
$$
w=s_{i_1}\ldots s_{i_k}\ldots s_{i_{\ell(w)}},\  i_1,\ldots i_{\ell(w)}\in I.
$$
It is known that the component of the differential $\D\til{M}_{\F}(w_1\cdot \lambda)\map
\D\til{M}_{\F}(w_2\cdot \lambda)$ is nontrivial if and only if $\ell(w_2)=\ell(w_1)+1$ and
$w_2$ follows $w_1$ in the Bruhat order on the Weyl group, i.e. there exists $k\le \ell(w_2)$ such that
$$
w_2=s_{i_1}\ldots s_{i_k}\ldots s_{i_{\ell(w_2)}}\text{ and }
w_1=s_{i_1}\ldots \hat s_{i_k}\ldots s_{i_{\ell(w_1)}}.
$$
Here $s_{i_1}\ldots \hat s_{i_k}\ldots s_{i_{\ell(w_1)}}$ denotes the expression for $w_2$ via the 
generators with the factor $s_{i_k}$ missing.

We start investigating components of the differentials in the complex $K_{\F}\bul(\lambda)$ 
with the simplest case.

\subsection{Components of the differential $d_i: \D\til{M}_{\F}(s_iw_0\cdot\lambda)\map
\D\til{M}_{\F}(w_0\cdot\lambda)$}
Below we provide an algebraic description of the components $d_i$ of the differential in the 
Grothendieck-Cousin complex on $\B_{\F}$ starting from the differential $d_{{\mathfrak {sl}}(2)}$
from the example above.

Denote by $\p_{i}$ (resp. by  $P_{i,\F}$) the standard $i$-th parabolic subalgebra with the unique 
negative root (resp. the corresponding subgroup in $G_{\F}$ with $\on{Lie}(P_{i,\F})=\p_{i}$). 
\vskip 1 mm
\noindent
\Rem
We formulate the next statement for the weight $0$. However a similar assertion remains true 
for any regular dominant integral weight $\lambda$. We restrict ourselves to the case of the the zero 
weight just to simplify notations in the proof.

\begin{thm} \label{square1}
\vskip 0.25 mm
\noindent
\vskip 0.25 mm\noindent\begin{itemize}
\item[(i)]
There exists a commutative square of $\U_{\F}(\g)$-modules as follows:
$$
\begin{array}{ccc}

\U_{\F}(\g)\ten_{\U_{\F}(\p_i)}\D\til{M}_{\F}^{{\mathfrak {sl}}(2)}(0)&\map&\D\til{M}_{\F}
(s_iw_0\cdot0)\\
\ &\ & \\
\downarrow{\lefteqn{\scriptstyle{\on{Ind}(d_{{\mathfrak {sl}}(2)})}}}&&\downarrow{\lefteqn{
\scriptstyle{d_i}}}\\
\ &\ & \\
\U_{\F}(\g)\ten_{\U_{\F}(\p_i)}\D\til{M}_{\F}^{{\mathfrak {sl}}(2)}(-2\rho)&\map&\D\til{M}_{\F}
(w_0\cdot0).
\end{array}
$$
\item[(ii)] The horizontal arrows in the square are isomorphisms.
\item[(iii)] The vertical arrows in the square are surjective.
\end{itemize}
\end{thm}

\begin{proof}
We prove the Theorem providing an explicit description of the commutative square in terms of local 
cohomology.

Note that in the $SL(2)$ case $\B_{\F}={\mathbb P}^1_{\F}$. Consider the diagonal action 
of the group $P_{i,\F}$ on $G_{\F}\times {\mathbb P}^1_{\F}$. Here $P_{i,\F}$ acts on the first
factor of the product via the right translations. It acts on the second factor via the projection
on its Levi factor $L_{i,\F}=SL(2,\F)$.

\begin{prop}
We have $G_{\F}\times {\mathbb P}^1_{\F}/P_{i,\F}\til\map\B_{\F}$.
\qed
\end{prop}

Consider the local cohomology spaces
$$
H^{\dim G+1}_{\{e\}\times pt}(G_{\F}\times{\mathbb P}^1_{\F},{\mathcal O}_{G_{\F}
\times{\mathbb P}^1_{\F}}) \text{ and }
H^{\dim G}_{\{e\}\times{\mathbb A}^1_{\F}}(G_{\F}\times{\mathbb P}^1_{\F},
{\mathcal O}_{G_{\F}
\times{\mathbb P}^1_{\F}}).
$$ 
\begin{lem}
We have
\vskip 0.25 mm\noindent\begin{itemize}
\item[(i)]
$\left(H^{\dim G+1}_{\{e\}\times pt}(G_{\F}\times{\mathbb P}^1_{\F},{\mathcal O}_{G_{\F}
\times{\mathbb P}^1_{\F}})\right)_{\U_{\F}(\p_i)\on{-coinv}}\til\map\U_{\F}(\g)\ten_{\U_{\F}(\p_i)}
\D\til{M}_{\F}^{{\mathfrak {sl}}(2)}(-2\rho);$
\item[(ii)]
$\left(H^{\dim G}_{\{e\}\times{\mathbb A}^1_{\F}}(G_{\F}\times{\mathbb P}^1_{\F},
{\mathcal O}_{G_{\F}
\times{\mathbb P}^1_{\F}})\right)_{\U_{\F}(\p_i)\on{-coinv}}\til\map\U_{\F}(\g)\ten_{\U_{\F}(\p_i)}
\D\til{M}_{\F}^{{\mathfrak {sl}}(2)}(0).$\qed
\end{itemize}
\end{lem}

Consider now the projection $m:\ G_{\F}\times {\mathbb P}^1_{\F}\map \B_{\F}$. 
\begin{lem}
There exist natural maps of $\U_{\F}(\g)$-modules
\vskip 0.25 mm\noindent\begin{itemize}
\item[(i)]
$
H^{\dim G+1}_{\{e\}\times pt}(G_{\F}\times{\mathbb P}^1_{\F},{\mathcal O}_{G_{\F}
\times{\mathbb P}^1_{\F}})\map H^{\dim \B}_{\{e\}}(\B_{\F},{\mathcal O}_{\B_{\F}});
$
\item[(ii)]
$H^{\dim G}_{\{e\}\times{\mathbb A}^1_{\F}}(G_{\F}\times{\mathbb P}^1_{\F},
{\mathcal O}_{G_{\F}
\times{\mathbb P}^1_{\F}})\map H^{\dim \B-1}_{C_{s_i,\F}}(\B_{\F},{\mathcal O}_{\B_{\F}}).
$
\end{itemize}
\end{lem}
\begin{proof}
Both statements of the lemma follow from Corollary~\ref{action}.
\end{proof}

\begin{lem}The maps from the above Lemma provide the isomorphisms:
\begin{gather*}
\left(H^{\dim G+1}_{\{e\}\times pt}(G_{\F}\times{\mathbb P}^1_{\F},{\mathcal O}_{G_{\F}
\times{\mathbb P}^1_{\F}})\right)_{\U_{\F}(\p_i)\on{-coinv}}\til\map 
H^{\dim \B}_{\{e\}}(\B_{\F},{\mathcal O}_{\B_{\F}}),\\
\left(H^{\dim G}_{\{e\}\times{\mathbb A}^1_{\F}}(G_{\F}\times{\mathbb P}^1_{\F},
{\mathcal O}_{G_{\F}
\times{\mathbb P}^1_{\F}})\right)_{\U_{\F}(\p_i)\on{-coinv}}\til\map
H^{\dim \B-1}_{C_{s_iw_0,\F}}(\B_{\F},{\mathcal O}_{\B_{\F}}).\qed
\end{gather*}
\end{lem}
Now all the four maps in the square in question are constructed, it remains to check that the square
 commutes. This follows from the functoriality of Grothendieck-Cousin complex 
(via the inclusion ${\mathbb P}^1_{\F}\til\map\overline{C}_{s_i,\F}\subset\B_{\F}$). The surjectivity 
of the vertical arrows follows from the one in the $SL(2)$ case and from the exactness of the induction 
functor.
\end{proof}

\subsection{The case of a simple reflection.}
Next we investigate the components in the differential of $K_{\F}\bul(\lambda)$ in the situation as follows.
Suppose we have a pair of elements of the Weyl group $w,w'\in W$ such that $\ell(w')=\ell(w)+1$
and $w'=ws_i$ for some $i\in I$. Consider the map $\D\til{M}(w\cdot\lambda)\map
\D\til{M}(w'\cdot\lambda)$. Again for the sake of simplicity we restrict ourselves to the case 
of the zero 
weight $\lambda$.
\begin{thm} \label{square2}
\vskip 0.25 mm
\noindent
\vskip 0.25 mm\noindent\begin{itemize}
\item[(i)]
There exists a commutative square of $\U_{\F}(\g)$-modules as follows:
$$
\begin{array}{ccc}

\Theta_{w,\F}\left(\U_{\F}(\g)
\ten_{\U_{\F}(\p_i)}\D\til{M}_{\F}^{{\mathfrak {sl}}(2)}(0)\right)&\map&\D\til{M}_{\F}
(w'w_0\cdot0)\\
\ &\ & \\
\downarrow{\lefteqn{\scriptstyle{\on{\Theta_{w,\F}\circ Ind}
(d_{{\mathfrak {sl}}(2)})}}}&&\downarrow{\lefteqn{
\scriptstyle{d_w'}}}\\
\ &\ & \\
\Theta_{w,\F}\left(\U_{\F}(\g)
\ten_{\U_{\F}(\p_i)}\D\til{M}_{\F}^{{\mathfrak {sl}}(2)}(-2\rho)\right)&\map&\D\til{M}_{\F}
(ww_0\cdot0).
\end{array}
$$
\item[(ii)] The horizontal arrows in the square are isomorphisms.
\item[(iii)] The vertical arrows in the square are surjective.

\end{itemize}
\end{thm}

\begin{proof}
We prove the Theorem providing an explicit description of the commutative square in terms of local 
cohomology. The proof goes along the lines of the one of Theorem~\ref{square1}.

Consider the action of the group $N_{w,\F}$  on $\B_{\F}$ and in particular the orbit
of the one dimensional Schubert cell $C_{s_iw_0,\F}$ under the action. 
\begin{lem}
When twisted by the action of $w\in W$ the set $N_{w,\F}\cdot C_{s_iw_0,\F}$ coincides with
the Schubert cell $C_{ws_iw_0,\F}$.
\qed
\end{lem}

 Again, like in the proof of Theorem~\ref{square1},  consider the diagonal action 
of the group $P_{i,\F}$ on $G_{\F}\times {\mathbb P}^1_{\F}$.
Consider the local cohomology spaces
$$
H^{\codim N_w}_{N_{w,\F}\times pt}(G_{\F}\times{\mathbb P}^1_{\F},{\mathcal O}_{G_{\F}
\times{\mathbb P}^1_{\F}}) \text{ and }
H^{\dim G}_{N_{w,\F}\times{\mathbb A}^1_{\F}}(G_{\F}\times{\mathbb P}^1_{\F},
{\mathcal O}_{G_{\F}
\times{\mathbb P}^1_{\F}}).
$$ 
\begin{lem}
We have
\begin{gather*}
T_{w,\F}\left(\left(H^{\dim G- \dim N_w+1}_{N_{w,\F}\times pt}
(G_{\F}\times{\mathbb P}^1_{\F},{\mathcal O}_{G_{\F}\times{\mathbb P}^1_{\F}})
\right)_{\U_{\F}(\p_i)\on{-coinv}}\right)\\
\til\map\Theta_{w,\F}\left(\U_{\F}(\g)\ten_{\U_{\F}(\p_i)}
\D\til{M}_{\F}^{{\mathfrak {sl}}(2)}(-2\rho)\right);\\
T_{w,\F}\left(\left(H^{\dim G- \dim N_w}_{N_{w,\F}\times{\mathbb A}^1_{\F}}
(G_{\F}\times{\mathbb P}^1_{\F},
{\mathcal O}_{G_{\F}
\times{\mathbb P}^1_{\F}})\right)_{\U_{\F}(\p_i)\on{-coinv}}\right)\\
\til\map\Theta_{w,\F}\left(\U_{\F}(\g)\ten_{\U_{\F}(\p_i)}
\D\til{M}_{\F}^{{\mathfrak {sl}}(2)}(0)\right).\qed
\end{gather*}
\end{lem}

Consider now the projection $m:\ G_{\F}\times {\mathbb P}^1_{\F}\map \B_{\F}$. 
\begin{lem}
There exist natural maps of $\U_{\F}(\g)$-modules
\vskip 0.25 mm\noindent\begin{itemize}
\item[(i)]
$
H^{\codim N_w+1}_{N_{w,\F}\times pt}(G_{\F}\times{\mathbb P}^1_{\F},{\mathcal O}_{G_{\F}
\times{\mathbb P}^1_{\F}})\map H^{\dim \B-\dim N_w}_{N_{w,\F}}(\B_{\F},
{\mathcal O}_{\B_{\F}});
$
\item[(ii)]
$H^{\codim N_w}_{N_{w,\F}\times{\mathbb A}^1_{\F}}(G_{\F}\times{\mathbb P}^1_{\F},
{\mathcal O}_{G_{\F}
\times{\mathbb P}^1_{\F}})\map H^{\dim \B- \dim N_w-1}_{N_{w,\F}\cdot C_{s_i,\F}}
(\B_{\F},{\mathcal O}_{\B_{\F}}).
$
\end{itemize}
\end{lem}
\begin{proof}
Both statements of the lemma follow from Corollary~\ref{action}.
\end{proof}
\begin{lem}The maps from the above Lemma provide the isomorphisms:
\begin{gather*}
\left(H^{\dim G-\dim N_w+1}_{N_{w,\F}\times pt}(G_{\F}\times{\mathbb P}^1_{\F},
{\mathcal O}_{G_{\F}
\times{\mathbb P}^1_{\F}})\right)_{\U_{\F}(\p_i)\on{-coinv}} \til\map 
H^{\dim \B-\dim N_w}_{N_{w,\F}\cdot\{e\}}(\B_{\F},{\mathcal O}_{\B_{\F}}),\\
\left(H^{\dim G-\dim N_w}_{N_{w,\F}\times{\mathbb A}^1_{\F}}(G_{\F}
\times{\mathbb P}^1_{\F},
{\mathcal O}_{G_{\F}
\times{\mathbb P}^1_{\F}})\right)_{\U_{\F}(\p_i)\on{-coinv}} \til\map
H^{\dim \B-\dim N_w-1}_{N_{w,\F}\cdot C_{s_iw_0,\F}}(\B_{\F},{\mathcal O}_{\B_{\F}}).
\end{gather*}
\qed
\end{lem}
Now applying the functor $T_{w,\F}$ to the isomorphisms from the above Lemma we
obtain isomorphisms as follows:
\begin{gather*}
\Theta_{w,\F}\left(\U_{\F}(\g)
\ten_{\U_{\F}(\p_i)}\D\til{M}_{\F}^{{\mathfrak {sl}}(2)}(0)\right)\til\map\D\til{M}_{\F}
(w'w_0\cdot0),\\
\Theta_{w,\F}\left(\U_{\F}(\g)
\ten_{\U_{\F}(\p_i)}\D\til{M}_{\F}^{{\mathfrak {sl}}(2)}(-2\rho)\right)\til\map\D\til{M}_{\F}
(ww_0\cdot0).
\end{gather*}
It remains to check that the square in question commutes. We use the same arguments as in the proof
of Theorem~\ref{square1}
\end{proof}
For any $w\in W$ consider the minimal length decomposition via simple reflections
$$w=s_{i_1}\ldots s_{i_{\ell(w)}},\ i_1,\ldots i_{\ell(w)}\in I.
$$
We have the sequence of the elements of the Weyl group
$$
v_0:=w,\ v_1=s_{i_1}\ldots s_{i_{\ell(w)-1}},\ldots
v_k=s_{i_1}\ldots s_{i_{\ell(w)-k}},\ldots,\ v_{\ell(w)}=e.
$$
The following statement is a consequence of the immediate generalization of Theorem~\ref{square2}
to the case of an arbitrary regular dominant integral weight $\lambda$.
\begin{cor} \label{chains}
We have the chains of the surjective maps as follows
\begin{gather*}
p_{e}=d_{s_{i_1}}:\ \D\til{M}_{\F}(\lambda)\map\D\til{M}_{\F}(s_{i_1}\cdot\lambda); \ldots\\
p_{v_{k}}=d_{v_k}\circ\ldots\circ d_{s_{i_1}}:\ \D\til{M}_{\F}(\lambda)\map
\D\til{M}_{\F}(s_{i_1}\cdot\lambda)\map\ldots\map
\D\til{M}(v_k\cdot\lambda);\\ \ldots;\ 
p_{v_0}=d_w\circ\ldots\circ d_{s_{i_1}}:\ \D\til{M}_{\F}(\lambda)\map\ldots\map
\D\til{M}_{\F}(w\cdot\lambda).
\end{gather*}
\end{cor}
\begin{proof}
Note that any two neighboring elements $v_k$ and $v_{k-1}$ differ by right multiplication by a 
simple reflection and their lengths differ by 1. Now apply Theorem~\ref{square2} (iii).
\end{proof}

\subsection{The case of a general component of the differential on $K_{\F}\bul(\lambda)$.}
Below we fix a reduced expression $w_0=s_{i_1}\ldots s_{i_{\ell(w_0)}}$ of the longest element in 
the 
Weyl group via the simple reflections. It is known that this provides reduced expressions for all the 
elements of the Weyl group.
\begin{lem}
Consider the complex contragradient dual for $K_{\F}\bul(\lambda)$. For any  $w\in W$ the maps 
$i_{w}=p_{v_0(w)}^*$ from  Lemma~\ref{chains} provide the embeddings 
$$\til{M}_{\F}(w\cdot\lambda)\map\til{M}_{\F}(\lambda)\ (=M_{\F}(\lambda))
$$
that do not depend on the chosen reduced expression of $w_0$.
\end{lem}
\begin{proof}
The statement of the Lemma follows from the general combinatorics of the Weyl group.
\end{proof}

Our aim now is to characterize the components of the differential
$\D\til{M}_{\F}(w\cdot\lambda)\map\D\til{M}_{\F}(w'\cdot\lambda)$, where
$\ell(w')=\ell(w)+1$ and $w'$ follows $w$ in the Bruhat order on the Weyl group. However we prefer 
to characterize the contragradient dual maps
$$\til{M}_{\F}(w'\cdot\lambda)\map\til{M}_{\F}(w\cdot\lambda).$$

\begin{thm} \label{embeddings}
\vskip 0.25mm
\noindent
\begin{itemize}
\item[(i)]
For any pair of elements $w,w'$ of the Weyl group such that $w'$ follows $w$ in the Bruhat order
the submodule $i_{w'}(\til{M}_{\F}(w'\cdot\lambda))\subset\til{M}_{\F}(\lambda)$ is embedded
into the submodule $i_w(\til{M}_{\F}(w\cdot\lambda))\subset\til{M}_{\F}(\lambda)$.
\item[(ii)] 
If in addition $\ell(w')=\ell(w)+1$ then the above embedding is provided by the component 
of the differential $\til{M}_{\F}(w'\cdot\lambda)\map\til{M}_{\F}(w\cdot\lambda)$
in the complex $\D K_{\F}\bul(\lambda)$.
\end{itemize}
\end{thm}
\begin{proof}
We begin the proof of the Theorem with generalizing the complex $\D K_{\F}\bul(\lambda)$. In fact 
the corresponding statement is contained in section 12 of \cite{Kempf}.

Recall that the Flag variety $G/B$ as well as the Schubert cells are defined over $\BZ$. Denote the 
corresponding objects by $\B_{\BZ}$ and $\{C_{w,\BZ}|w\in W\}$. The standard line bundles 
$\CL(\lambda)$ are also defined over $\BZ$.

The following statement is due to Kempf.
\begin{prop} \label{Kempf}
\vskip 0.25 mm
\noindent
\begin{itemize}
\item[(i)]
There exists a complex of free $\BZ$-modules $K_{\BZ}\bul(\lambda)$ with
$$
K_{\BZ}^m(\lambda)=\underset{\ell(w)=m}\bigoplus\ H_{C_{ww_0,\BZ}}^m(\B_{\BZ},
\CL(\lambda)).
$$
\item[(ii)] Higher cohomology groups of the complex vanish and $H^0(K_{\BZ}\bul(\lambda))=
H^0(\B_{\BZ},\CL(\lambda))$.

\item[(iii)] For any finite field $\f$ we have 
$K_{\BZ}\bul(\lambda)\ten_{\BZ}\F=K_{\F}\bul(\lambda)$.\qed
\end{itemize}
\end{prop}
\vskip 1mm
\noindent
\Rem
The complex from the above Proposition is called the global Grothendieck-Cousin complex
for $\CL(\lambda)$ on $\B_{\BZ}$.

Note also that  $K_{\BZ}\bul(\lambda)\ten_{\BZ}\BC$ equals the 
contragradient Bernstein-Gelfand-Gelfand complex for the regular dominant integral weight 
$\lambda$.

We prefer to work with the contragradient dual complex
$\D K_{\BZ}\bul(\lambda)$. The duality makes sense since $K_{\BZ}\bul(\lambda)$ consists of
$\BZ$-free modules graded by the torus action so that the grading components are of finite rank 
over $\BZ$.

For a given element of the Weyl group $w$ (and the fixed reduced expression of $w$ via simple 
reflections $w=s_{i_1}\ldots s_{i_{\ell(w)}}$) we have the  chain of maps 
$$
\D H_{C_{ww_0,\BZ}}^{\ell(w)}(\B_{\BZ},
\CL(\lambda))\map\D H_{C_{ws_{i_{\ell(w)}}w_0,\BZ}}^{\ell(w)}(\B_{\BZ},
\CL(\lambda))\map\ldots\map \D H_{C_{w_0,\BZ}}^{\ell(w)}(\B_{\BZ},
\CL(\lambda))
$$
given by the components of the differentials in $\D K_{\BZ}\bul(\lambda)$.
\begin{lem}
All the arrows  in the chain are embeddings.
\end{lem}
\begin{proof}
The corresponding statements over all the  fields $\F$ are proved above. It follows that the 
assertion of the Lemma holds since the property $$\{\text{\sf {not to be an embedding}}\}$$ 
is a closed one.
\end{proof}
\begin{cor}
The $\BZ$-modules $\on{Im}\left(\D H_{C_{ww_0,\BZ}}^{\ell(w)}(\B_{\BZ},
\CL(\lambda))\right)\subset \D H_{C_{w_0,\BZ}}^{\ell(w)}(\B_{\BZ},
\CL(\lambda))$ do not depend on the choice of the reduced expression for $w_0$.
\end{cor}
\begin{proof}
This also follows from the corresponding statement for the finite fields.
\end{proof}
Let us now sum up what we have got already. We have a number of  families  over $\on{Spec Z}$ 
of submodules
in $\D H_{C_{w_0,\BZ}}^{\ell(w)}(\B_{\BZ},
\CL(\lambda))$ that is viewed as a quasicoherent sheaf over  $\on{Spec Z}$.
Each family is enumerated  by an element  of the Weyl group. We know also that over the generic 
point of $\on{Spec Z}$ the fiber of the family enumerated by $w'$ is a submodule in the fiber 
of the family enumerated by $w$ if and only if $w'$ follows $w$ in the Bruhat order. This is a 
well-known property of the Bernstein-Gelfand-Gelfand resolution. But the property
$$\{\text{\sf {to be a submodule}}\}$$
is a closed property in a flat family. Thus we obtain the statement as follows.

\begin{prop}
For any finite field $\f$ we have and a pair of elements $w,w'$  of the Weyl group such that $w'$ 
follows $w$ in the Bruhat order we have
$$\on{Im}(\til{M}_{\F}(w'\cdot\lambda))\hookrightarrow\on{Im}(\til{M}_{\F}(w\cdot\lambda))
$$
as submodules of $\til{M}_{\F}(\lambda)$. \qed
\end{prop}
This essentially finishes the proof of the first part of the Theorem.We leave the proof of the second 
part to the reader.
\end{proof}
\vskip 1mm
\noindent
\Rem
Note that the second part of the above Theorem provides in particular an implicit characterization
of the components of the differential in the complex
$K_{\F}\bul(\lambda)$: 
$$ \D\til{M}_{\F}(w\cdot\lambda)\map\D\til{M}_{\F}
(w'\cdot\lambda)\text{ for }\ell(w')=\ell(w)+1.
$$
Namely these components can be deduced from the structure of the lattice of 
$\U_{\F}(\g)$-submodules of $\til{M}(\lambda)$. This lattice of submodules is obtained using 
the properties of the special cases of the differential (see Theorem~\ref{square1} and 
Theorem~\ref{square2}).

\section{Comparing the quasi-BGG complex with $K_{\F}\bul(\lambda)$.}
We begin this section with recalling the quantum group setting and presenting briefly the construction
of the quasi-BGG complex from \cite{Ar0}.

Fix a {\em Cartan
datum} $(I,\cdot)$ of the finite type and a {\em simply connected root
datum} $(Y,X,\ldots)$ of the type $(I,\cdot)$.
Below we denote the {\em Drinfeld-Jimbo
quantum group} defined over the field $\Q(v)$ of rational functions in
$v$ (resp.  the {\em Lusztig version of the quantum group} defined
over $\CA=\BZ[v,v^{-1}]$) that correspond to the root data by  $\U$ (resp.  by $\U_{\CA}$).

Fix the natural triangular decompositions of the algebra $\U_{\CA}$
 as follows:
$\U_{\CA}=\U_{\CA}^-\ten\U_{\CA}^0\ten\U_{\CA}^+$, where the
positive (resp.  negative) subalgebras are generated by the quantum
divided powers of the positive (resp.  negative) root generators in
the corresponding algebra.  We call the subalgebra
$\U_{\CA}^+\ten\U_{\CA}^0$  the
{\em positive quantum Borel subalgebra} in $\U_{\CA}$ 
 and denote it by $\BB_{\CA}^+$.  The
negative Borel subalgebra $\BB_{\CA}^-$ is defined in
a similar way.

\subsection{Twisted quantum parabolic subalgebras in $\U_\CA$} Recall
that Lusztig has constructed an action of the {\em braid group}
$\mathfrak B$ corresponding to the Cartan data $(I,\cdot)$ by
automorphisms of the quantum group $\U_{\CA}$ well defined with
respect to the $X$-gradings (see \cite{L1}, Theorem 3.2).  Fix a
reduced expression of the maximal length element $w_0\in W$ via the
simple reflection elements:  $$ w_0=s_{i_1}\ldots
s_{i_{\sharp(R^+)}},\ i_k\in I.  $$ It has been mentioned already in the previous section
 that this reduced
expression provides reduced expressions for all the elements $w\in W$:
$ w=s_{i_1^w}\ldots s_{i_{l(w)}^w},\ i_k\in I.  $

Consider the standard generators $\{T_i\}_{i\in I}$ in the braid group
$\mathfrak B$.  Lifting the reduced expressions for the elements $w$
from $W$ into $\mathfrak B$ we obtain the set of elements in the braid
group of the form $T_w:= T_{i_1^w}\ldots T_{i_{l(w)}^w}$.

In particular we obtain the set of {\em twisted Borel subalgebras}
$w(\BB_\CA^+)=T_w(\BB_\CA^+)\subset\U_\CA$.  Note that
$w_0(\BB_\CA^+)=\BB_\CA^-=\U_\CA^-\ten\U_\CA^0$.

Fix a subset $J\subset I$ and consider the {\em quantum parabolic
subalgebra} $\P_{J,\CA}\subset\U_\CA$.  By definition this subalgebra
in $\U_\CA$ is generated over $\U_\CA^0$ by the elements $E_i$, $i\in
I$, $F_j$, $j\in J$, and by their quantum divided powers.  The
previous construction provides the set of {\em twisted quantum
parabolic subalgebras} $w(\P_{J,\CA}):=T_w(\P_{J,\CA})$ of the type
$J$ with the twists $w\in W$.

Note that the triangular decomposition of the algebra $\U_\CA$
provides the ones for the algebras $w(\BB_\CA^+)$ and $w(\P_{J,\CA})$:
$$ w(\BB_\CA^+) = (w(\BB_\CA^+))^- \ten \U_\CA^0 \ten (w(\BB_\CA^+))^+
\text{ and } w(\P_{J,\CA}) = (w(\P_{J,\CA}))^- \ten \U_\CA^0 \ten
(w(\P_{J,\CA})^+, $$ where $\left(w(\BB_\CA^+)\right)^+ =
w(\BB_\CA^+)\cap \U_\CA^-$, $\left(w(\P_{J,\CA})\right)^+ =
w(\P_{J,\CA})\cap \U_\CA^-$ etc.

Fix a dominant integral weight $\lambda\in X$.  Consider the
module over the  quantum group $\U_{\CA}$ given by $$ \D
W_{\CA}(\lambda):=\left(\Coind_{\BB^-_{\CA}}^{\U_{\CA}}\BC(\lambda)
\right)^{\operatorname{fin}} \text{ (resp.  by }
W_{\CA}(\lambda):=\left(\Ind_{\BB^+{\CA}}^{\U_{\CA}}\BC(\lambda)
\right)_{\operatorname{fin }}).  $$ Here $(*)^{\operatorname{fin}}$
(resp.  $(*)_{\operatorname{fin}}$) denotes the maximal finite
dimensional submodule (resp.  quotient module) in $(*)$.  The module
$\D W_{\CA}(\lambda)$ (resp.  $W_{\CA}(\lambda)$) is called {\em the
contragradient Weyl module} (resp.  {\em the Weyl module}) over
$\U_{\CA}$ with the highest weight $\lambda$.

\subsection{Semiinfinite induction and coinduction} From now on we
will use freely the technique of associative algebra semiinfinite
homology and cohomology for a graded associative algebra $A$ with two
subalgebras $B,N\subset A$ equipped with a triangular decomposition
$A=B\ten N$ on the level of graded vector spaces.  We will not recall
the construction of these functors referring the reader to \cite{Ar1}
and \cite{Ar2}.

Let us mention only that these functors are bifunctors ${\mathsf
D}(A\mod)\times{\mathsf D}(\oppA\mod)\map{\mathsf D}({\mathcal V}ect)$
where the associative algebra $\oppA$ is defined as follows.

Consider the semiregular $A$-module $S_A^N:=A\ten_NN^*$.  It is proved
in \cite{Ar2} that under very weak conditions on the algebra $A$ the
module $S_A^N$ is isomorphic to the $A$-module
$\left(S_A^N\right)':=\hom_B(A,B)$.  Thus $\End_A(S_A^N)\supset
N^\opp$ and $\End_A(S_A^N)\supset B^\opp$ as subalgebras.  The algebra
$\oppA$ is defined as the subalgebra in $\End_A(S_A^N)$ generated by
$B^\opp$ and $N^\opp$.  It is proved in \cite{Ar2} that the algebra
$\oppA$ has a triangular decomposition $\oppA=N^\opp\ten B^\opp$ on
the level of graded vector spaces.  Yet for an arbitrary algebra $A$
the algebras $\oppA$ and $A^\opp$ do not coincide.

However the following statement shows that in the case of quantum
groups that correspond to the root data $(Y,X,\ldots)$ of the {\em
finite} type $(I,\cdot)$  the equality of $A^\opp$ and $\oppA$ holds.

\begin{prop} (see Proposition 3.3.1 from \cite{Ar0})
We have \begin{itemize} 
\item[(i)] $\U^{\sharp}=\U^\opp$,
$\U_\CA^{\sharp}=\U_\CA^\opp$;
\item[(ii)] $w(\BB^+)^\sharp= w(\BB^+)^\opp$, $w(\BB_\CA^+)^\sharp=
w(\BB_\CA^+)^\opp$, 
$w(\P_{J,\CA})^\sharp= w(\P_{J,\CA})^\opp$. \qed
\end{itemize}
\end{prop}

\begin{defn}
Let $M\bul$ be a convex complex of $w(\BB_\CA^+)$-modules.
By definition set \begin{gather*}
\SInd_{w(\BB_\CA^+)}^{\U_\CA}(M\bul):=
\tor_{\si+0}^{w(\BB_\CA^+)}(S_{\U_\CA}^{\U_\CA^+},M\bul) \text{ and }\\
\SCoind_{w(\BB_\CA^+)}^{\U_\CA}(M\bul):=
\ext^{\si+0}_{w(\BB_\CA^+)}(S_{\U_\CA}^{\U_\CA^-},M\bul).
\end{gather*} 
\end{defn}
\begin{lem}
(see \cite{Ar4}) \begin{itemize} \item[(i)]
$\tor_{\si+k}^{w(\BB_\CA^+)}(S_{\U_\CA}^{\U_\CA^+},\cdot)=0$ for
$k\ne0$; \item[(ii)]
$\ext^{\si+k}_{w(\BB_\CA^+)}(S_{\U_\CA}^{\U_\CA^-},\cdot)=0$ for
$k\ne0$; \item[(iii)] $\SInd_{w(\BB_\CA^+)}^{\U_\CA}(\cdot)$ and
$\SCoind_{w(\BB_\CA^+)}^{\U_\CA}(\cdot)$ define {\em exact} functors
$$w(\BB_\CA^+)\mod\map\U_\CA\mod.\qed$$ \end{itemize}
\end{lem}

 Similar statements hold for the algebras $w(\P_{J,\CA})$.

\subsubsection{quasi-Verma modules over $\U_{\CA}$.}
We define the {\em quasi-Verma
module} over the algebra $\U_\CA$ with the highest
weight $w\cdot\lambda$ by $$ M_\CA^w(w\cdot \lambda):=
\SInd_{w(\BB_\CA^+)}^{\U_\CA}(\CA(\lambda)).  $$
The {\em
contragradient quasi-Verma module} $ \D M_\CA^w(w\cdot\lambda)$ is defined by 
$$ \D M_\CA^w(w\cdot\lambda):=
\SCoind_{w(\BB_\CA^+)}^{\U_\CA}(\CA(\lambda)).
$$ 
\vskip 1mm
\noindent
\Rem
Note that the definitition of a contragradient quasi-Verma module over $\U_{\CA}$ looks 
different from the one of the similar object over $\U_{\F}(\g)$. We will have to overcome this 
inconvenience later.

We list the main properties of quasi-Verma modules.

\begin{prop}
(see \cite{Ar4}) \begin{itemize} \item[(i)] Fix a dominant
integral weight $\lambda\in X$.  Suppose that $\xi\in\BC^*$ is not a
root of unity.  Then the $\U_\xi$-module
$M_\xi^w(w\cdot\lambda):=M_\CA^w(w\cdot\lambda)\ten_{\CA}\BC$ (resp.
$\D M_\xi^w(w\cdot\lambda):=\D M_\CA^w(w\cdot\lambda)\ten_{\CA}\BC$)
is isomorphic to the usual Verma module $M_\xi(w\cdot\lambda)$ (resp.
to the usual contragradient Verma module $\D M_\xi(w\cdot\lambda)$).
\item[(ii)] For any $\lambda\in X$ we have \begin{gather*}
\ch(M_\CA^w(w\cdot\lambda)) = \ch(\D M_\CA^w(w\cdot\lambda)
  =
\frac{e^{w\cdot\lambda}}{\prod_{\alpha\in R^+}(1-e^{-\alpha})}.\qed
\end{gather*} \end{itemize}
\end{prop}

 Thus for a dominant weight $\lambda$ one
can consider $M_\CA^w(w\cdot\lambda)$ as a flat family of modules over
the quantum group for various values of the quantizing parameter with
the fiber at the generic point equal to the Verma module
$M_\xi(w\cdot\lambda)$.

\subsubsection{Quasi-BGG complex in the case of $\U_{\CA}({\mathfrak {sl}}(2))$.}
In \cite{Ar0} this case was investigated throughly and the following statement was proved.

\begin{lem} \label{sl2case}
For every positive integer $\mu$ there exists an exact complex of $\U_{\CA}
({\mathfrak {sl}}_2)$-modules
$$ 0\map
M_\CA^s(s\cdot \mu)\map M_\CA^e(\mu)\map W_\CA(\mu)\map0.\qed $$
\end{lem}
We call the complex $M_\CA^s(s\cdot \mu)\map M_\CA^e(\mu)$ the quasi-BGG 
complex for the weight $\mu$ and denote it by $B_\CA\bul(\mu)$.

\subsection{Construction of the quasi-BGG complex for general $\U_{\CA}$.} 
 Here we extend
the previous considerations to the case of the quantum group $\U_\CA$
for arbitrary root data $(Y,X,\ldots)$ of the finite type $(I,\cdot)$.
Fix a dominant weight $\mu \in X$.

First we construct an inclusion
$M_\CA^{w'}(w'\cdot\mu)\hookrightarrow
M_\CA^{w}(w\cdot\mu)$ for a pair of elements $w',w\in W$ such
that $\ell(w')=\ell(w)+1$ and $w'$ follows $w$ in the Bruhat order on the Weyl group.
In fact we can do it explicitly only for $w'$ and $w$ differing by a
simple reflection:  $w'=ws_i$, $i\in I$.

Consider the twisted quantum parabolic subalgebra $w(\P_{i,\CA})$.
Then $w(\P_{i,\CA})\supset w(\BB_\CA^+)$ and $w(\P_{i,\CA})\supset
ws_i(\BB_\CA^+)$.  Consider also the Levi quotient algebra
$w(\P_{i,\CA})\map w(\LL_{i,\CA})$.  The algebra $w(\LL_{i,\CA})$
is isomorphic to $\U_\CA({\mathfrak {sl}}_2)\ten_{\U_\CA^0({\mathfrak
{sl}}_2)}\U_\CA^0$.

By Lemma \ref{sl2case} we have a natural inclusion of $w
(\LL_{i,\CA})$-modules
$$\SInd_{ws_i(\LL_{i,\CA}^+)}^{w(\LL_{i,\CA})}\CA(\mu)
\hookrightarrow
\SInd_{w(\LL_{i,\CA}^+)}^{w(\LL_{i,\CA})}\CA(\mu).$$ 
\begin{lem}
\vskip 0.25mm
\noindent
\begin{itemize} 
\item[(i)] $\SInd_{w(\BB_\CA^+)}^{\U_\CA}(\cdot)=
\SInd_{w(\P_{i,\CA})}^{\U_\CA}\circ
\SInd_{w(\BB_\CA^+)}^{w(\P_{i,\CA})}(\cdot)$.  
\item[(ii)]
$\SInd_{ws_i(\BB_\CA^+)}^{\U_\CA}(\cdot)=
\SInd_{w(\P_{i,\CA})}^{\U_\CA}\circ
\SInd_{ws_i(\BB_\CA^+)}^{w(\P_{i,\CA})}(\cdot)$.  
\item[(iii)]
$\SInd_{w(\BB_\CA^+)}^{\U_\CA}(\CA(\mu))=
\SInd_{w(\P_{i,\CA})}^{\U_\CA}\circ
\Res_{w(\P_{i,\CA})}^{w(\LL_{i,\CA})}\circ
\SInd_{W(\LL_{i,\CA}^+)}^{w(\LL_{i,\CA})}(\mu)$.
\item[(iv)] 
$\SInd_{ws_i(\BB_\CA^+)}^{\U_\CA}(\CA(\mu))=
\SInd_{w(\P_{i,\CA})}^{\U_\CA}\circ
\Res_{w(\P_{i,\CA})}^{w(\LL_{i,\CA})}\circ
\SInd_{ws_i(\LL_{i,\CA}^+)}^{w(\LL_{i,\CA})}(\mu)$.  \qed
\end{itemize}
\end{lem}
 \begin{cor} For $w'=ws_i>w$ in the Bruhat order we
have a natural inclusion of $\U_\CA$-modules $i_\CA^{ws_i,w}:\
M_\CA^{ws_i}(ws_i\cdot\mu)\hookrightarrow
M_\CA^w(w\cdot\mu).  \qed $
\end{cor}

Recall that if $v$ acts on $\BC$ by 
$\xi$ that is  not a root of unity then the $\U_\xi$-module 
$M_\xi^w(w\cdot\mu):=M_\CA^w(w\cdot\mu)\ten_\CA\BC$ is
isomorphic to the usual Verma module $M_\xi(w\cdot\mu)$.  Thus
the morphism $i_\xi^{w',w}$ coincides with the standard inclusion of
Verma modules constructed by J.  Bernstein, I.M.  Gelfand and S.I.
Gelfand in \cite{BGG} that becomes a component of the differential in
the BGG resolution.  In other words we see that the flat family of
inclusions $i_\xi^{ws_i,w}:\
M_\ell^{ws_i}(ws_i\cdot\mu)\hookrightarrow
M_\xi^w(w\cdot\mu) $ defined   for  $\xi\in\BC^*\setminus\{$roots of
unity$\}$ can be extended naturally over the whole $\Spec \CA$.

Iterating the inclusion maps we obtain a flat family of submodules
$i^w_\xi(M_\xi^w(w\cdot\mu))\subset M_\xi^e(\mu)$ for
$\xi\in\Spec \CA$, $w\in W$, providing an extension of the standard
lattice of Verma submodules in $M_\xi(\mu)$ defined a priori for
$\xi\in\BC^*\setminus\{$roots of unity$\}$.

\begin{lem}  For a pair of elements $w',w\in W$ such that
$\ell(w')=\ell (w)+1$ and $w'$ follows $w$ in the Bruhat order we have $$
i_\CA^{w'}(M_\CA^{w'}(w'\cdot\mu)) \hookrightarrow
i_\CA^w(M_\CA^{w}(w\cdot\mu)).  \qed$$ 
\end{lem}

Now using the standard combinatorics of the classical BGG resolution
we obtain the following statement.

\begin{thm} (see Theorem 3.6.2 in \cite{Ar0})
 There exists a complex of $\U_\CA$-modules
$B\bul_\CA(\mu)$ with $$
B_\CA^{-k}(\mu)=\underset{w\in W,\ell(w)=k}{\bigoplus}
M_\CA^w(w\cdot\mu) $$ and with differentials provided by
direct sums of the inclusions $i_\CA^{w',w}$.  \qed
\end{thm}
\begin{defn} We call the complex $B\bul_\CA(\mu)$ the {\em
quasi-BGG complex} for the regular dominant integral weight $\mu\in X$.
\end{defn}

\subsection{Main theorem.}
Like in~\cite{L1} 
consider the specialization of $\U_\CA$ in characteristic $p$. 
Namely let $\CA_p'$ be the quotient of $\CA$ by the ideal generated by the
 $p$-th cyclotomic polynomial. Then 
$\CA_p'/(v-1)$ is isomorphic to the finite field ${\mathbb F}_p$. 
Thus the algebraic closure $\F$
becomes a $\CA$-algebra. We set $\U_{\F}:=\U_\CA\ten_{\CA}\F$. 
 It is known that
the quotient of the algebra $\U_{\F}$ by certain central elements is isomorphic to 
$\U_{\BZ}(\g)\ten \F$, where, as before, $\U_\BZ(\g)$ denotes the Kostant integral form 
for the universal enveloping algebra of 
$\g$. 

\begin{lem}
The algebra $\U_{\F}(\g)$ acts on the complex $B\bul_\CA(\lambda)\ten_{\CA}\F=:
B\bul_{\F}(\lambda)$.
\end{lem}
\begin{proof}
First one checks that for a regular dominant integral weight $\lambda$ the module
$M_\CA^e(\lambda)\ten_{\CA}\F$ is isomorphic to $M_{\F}(\lambda)$ and thus 
the algebra $\U_{\F}(\g)$ acts on it. Next note that all the modules 
$M_\CA^w(w\cdot\lambda)\ten_{\CA}\F$ can be viewed as submodules in
 $M_\CA^e(\lambda)\ten_{\CA}\F$, thus they are acted by $\U_{\F}(\g)$ as well. Finally the 
components of the differentials in $B\bul_\CA(\lambda)\ten_{\CA}\F$ are just the inclusions
$M_\CA^{w'}(w'\cdot\lambda)\ten_{\CA}\F\map M_\CA^w(w\cdot\lambda)\ten_{\CA}\F$
for $w,w'\in W$, such that $\ell(w')=\ell(w)+1$ and $w'$ follows $w$ in the Bruhat order.
\end{proof}

Now we are ready to formulate the main result of the whole paper.

\begin{thm} \label{main}
For a regular dominant integral weight $\lambda$ the complexes $\D B\bul_{\F}(\lambda)$
and $K_{\F}\bul(\lambda)$ are isomorphic as complexes of $\U_{\F}(\g)$-modules.
\end{thm}
\begin{proof}
Formally note that both complexes look as follows. As a $\U_{\F}(\g)$-module
each of them is a direct some of certain $\U_{\F}(\g)$-modules enumerated by the Weyl group, and
the characters of the corresponding modules in the direct sums coincide.

The following statement is the key one in the proof of the Theorem.
\begin{prop} \label{spec}
The specialization of the quasi-Verma module $M^w_{\CA}(w\cdot\lambda)$ into $\F$
is isomorphic to the $\U_{\F}(\g)$-quasi-Verma module $\til{M}_{\F}(w\cdot\lambda)$.
\end{prop}
\begin{proof}
Again we use the technique of semiinfinite cohomology with few comments and with references to 
\cite{Ar2}.
Note that the algebra  $\U_{\F}(\g)$ posesses a triangular decomposition
$\U_{\F}(\g)=\U_{\F}(\n^-)\ten\U_{\F}(\h)\ten\U_{\F}(\n^+)$ on the level of $\F$-vector spaces.
This decomposition comes as the specialization of the one of $\U_{\CA}$ into the field $\F$.

\begin{lem} \label{2}
We have
\begin{itemize}
\item[(i)] 
The semiregular modules 
$$S_{\U_{\F}(\g)}^{\U_{\F}(\n^+)}=S_{\U_{\CA}}^{\U_{\CA}^+}\ten_{\CA}
\F\text{ and }S_{\U_{\F}(\g)}^{\U_{\F}(\n^-)}=S_{\U_{\CA}}^{\U_{\CA}^-}\ten_{\CA}\F.$$
\item[(ii)] 
We have the isomorphism of algebras $\U_{\F}(\g)^\sharp\til\map\U_{\F}(\g)$.
\item[(iii)] 
For any $\CA$-free module $M$ we have an isomorphisms of the $\U_{\F}(\g)$-modules
\begin{gather*}
\SInd_{w(\BB_\CA^+)}^{\U_\CA}(M)\ten_{\CA}\F\til\map
\SInd_{\U_{\F}(w(\b^+))}^{\U_{\F}(\g)}(M\ten_{\CA}\F)\text{ and }\\
\SCoind_{w(\BB_\CA^+)}^{\U_\CA}(M)\ten_{\CA}\F\til\map
\SCoind_{\U_{\F}(w(\b^+))}^{\U_{\F}(\g)}(M\ten_{\CA}\F).\qed
\end{gather*}
\end{itemize}
\end{lem}

Thus we are to compare $\D\til{M}_{\F}(w\cdot\lambda)$ and 
$\SCoind_{\U_{\F}(w(\b^+))}^{\U_{\F}(\g)}(\F(\lambda))$. Let us spell out the definition of the last
module in down-to-earth terms.

\begin{lem} \label{1}
The $\U_{\F}(\g)$-module $\SCoind_{\U_{\F}(w(\b^+))}^{\U_{\F}(\g)}(\F(\lambda))$
is isomorphic to 
$$\hom_{\U_{\F}(w(\b^+))}(\F(-\lambda),(\U_{\F}(w(\b^+)\cap\n^-))^*
\ten_{\U_{\F}(w(\b^+)\cap\n^-)}\U_{\F}(\g)\ten_{\U_{\F}(\n^+)}(\U_{\F}(\n^+))^*).$$
\qed
\end{lem}

Here $(\cdot)^*$ denotes the restricted dual module for $(\cdot)$.

\begin{prop} \label{4}
There exists a natural map of $\U_{\F}(\g)$-modules
$$
\alpha:\ \D\til{M}(w\cdot\lambda)\map\SCoind_{\U_{\F}(w(\b^+))}^{\U_{\F}(\g)}(\F(\lambda)).
$$
\end{prop}
\begin{proof}Below we cheat a little forgetting about $T_{w,\F}$ in the definition of 
$\D\til{M}(w\cdot\lambda)$.
Note that we have a canonical identification
\begin{gather*}
\hom_{\U_{\F}(\g)}\left({\mathcal S}_{w,\F}\ten_{\U_{\F}(\g)}M_{\F}(w_0\cdot\lambda),
\SCoind_{\U_{\F}(w(\b^+))}^{\U_{\F}(\g)}(\F(\lambda))\right) \til\map\\
\hom_{\U_{\F}(\g)}\left(M_{\F}(w_0\cdot\lambda),\hom_{\U_{\F}(\g)}\left({\mathcal S}_{w,\F},
\SCoind_{\U_{\F}(w(\b^+))}^{\U_{\F}(\g)}(\F(\lambda))\right)\right) \til\map\\
\hom_{\U_{\F}(\g)}\left(M_{\F}(w_0\cdot\lambda),\hom_{\U_{\F}(\n_w^-)}\left((\U_{\F}(\n_w^-))^*,
\SCoind_{\U_{\F}(w(\b^+))}^{\U_{\F}(\g)}(\F(\lambda))\right)\right).
\end{gather*}
Now we use the following statement.
\begin{lem}
There exists a natural $\U_{\F}(\g)$-module map
$$
M_{\F}(w_0\cdot\lambda)\map \hom_{\U_{\F}(\n_w^-)}\left(\U_{\F}(\n_w^-))^*,
\SCoind_{\U_{\F}(w(\b^+))}^{\U_{\F}(\g)}(\F(\lambda))\right).
\qed
$$
\end{lem}
Thus we obtain a canonical element
\begin{gather*}
\on{Id}\in \hom_{\U_{\F}(\g)}(M_{\F}(w_0\cdot\lambda),M_{\F}(w_0\cdot\lambda))\\ \map
\hom_{\U_{\F}(\g)}\left({\mathcal S}_{w,\F}\ten_{\U_{\F}(\g)}M_{\F}(w_0\cdot\lambda),
\SCoind_{\U_{\F}(w(\b^+))}^{\U_{\F}(\g)}(\F(\lambda))\right).
\end{gather*}
This way we obtain the required map $\alpha$.
\end{proof}
\begin{prop} \label{3}
The map $\alpha$ is an isomorphism of $\U_{\F}(\g)$-modules.\qed
\end{prop}

Now the statement of Proposition~\ref{spec} follows from Lemma~\ref{2}, Lemma~\ref{1},
Proposition~\ref{4} and Proposition~\ref{3}.
\end{proof}

It remains to compare the differentials in the complexes $\D B\bul_{\F}(\lambda)$
and $K_{\F}\bul(\lambda)$. To do this note that the components of the differentials 
in the complexes $B\bul_{\F}(\lambda)$ and $\D K_{\F}\bul(\lambda)$ are obtained by 
the same algorythm. Namely we have two lattices of submodules in $M_{\F}(\lambda)$.
Both lattices are enumerated by the Weyl group and the inclusions of submodules in the lattices
agree with the Bruhat order on $W$. In both cases for $w,w'\in W$, such that $\ell(w')=
\ell(w)+1$ and $w'$ follows $w$ in the Bruhat order, we obtain the component of the differential in the
 corresponding complex by taking the inclusion map in the corresponding lattice.

Now the statement below follows form Theorem~\ref{embeddings} and from the construction 
of the quasi-BGG complex from the beginning of the present section.

\begin{lem}
The lattices of submodules $$\{\on{Im}(\til{M}_{\F}(w\cdot\lambda))\}_{w\in W}\text{ and }
\{\on{Im}(M^w_{\CA}(w\cdot\lambda)\ten_{\CA}\F\}_{w\in W}$$ in $M_{\F}(\lambda)$ coincide.
\qed
\end{lem}
The Lemma completes the proof of Theorem~\ref{main}.
\end{proof}

\end{document}